\documentclass[12pt]{article}
\usepackage{amssymb,amsmath,bm}
\usepackage{mathrsfs}
\usepackage{constants}
\usepackage{color,verbatim}
\usepackage{enumerate}
\usepackage{appendix}
\usepackage{authblk}
\usepackage[usenames,dvipsnames]{xcolor}

\usepackage{hyperref}
\hypersetup{
    colorlinks,
    citecolor=black,
    filecolor=black,
    linkcolor=black,
    urlcolor=black
}

\begin{document}
\newcommand{\bea}{\begin{eqnarray}}
\newcommand{\ena}{\end{eqnarray}}
\newcommand{\beas}{\begin{eqnarray*}}
\newcommand{\enas}{\end{eqnarray*}}
\newcommand{\beq}{\begin{equation}}
\newcommand{\enq}{\end{equation}}
\def\qed{\hfill \mbox{\rule{0.5em}{0.5em}}}
\newcommand{\bbox}{\hfill $\Box$}
\newcommand{\ignore}[1]{}
\newcommand{\ignorex}[1]{#1}
\newcommand{\wtilde}[1]{\widetilde{#1}}
\newcommand{\mq}[1]{\mbox{#1}\quad}
\newcommand{\bs}[1]{\boldsymbol{#1}}
\newcommand{\qmq}[1]{\quad\mbox{#1}\quad}
\newcommand{\qm}[1]{\quad\mbox{#1}}
\newcommand{\nn}{\nonumber}
\newcommand{\Bvert}{\left\vert\vphantom{\frac{1}{1}}\right.}
\newcommand{\To}{\rightarrow}
\newcommand{\supp}{\mbox{supp}}
\newcommand{\law}{{\cal L}}
\newcommand{\mS}{\mathcal{S}}
\newcommand{\Z}{\mathbb{Z}}
\newcommand{\jcolor}[1]{\textcolor{ForestGreen}{#1}}  % Jay's colored text
\newcommand{\jcomm}[1]{\marginpar{\tiny\jcolor{#1}}}  % Jay's marginpars
\newcommand{\ucolor}[1]{\textcolor{red}{#1}}  % Umit's colored text
\newcommand{\ucomm}[1]{\marginpar{\tiny\ucolor{#1}}}  % Umit's marginpars
\newcommand{\lcolor}[1]{\textcolor{blue}{#1}}  % Larry's colored text
\newcommand{\lcomm}[1]{\marginpar{\tiny\lcolor{#1}}}  % Larry's marginpars

\newcounter{monsterThmCounter} % for continuing the list in the monster theorems

%Erd\H{o}s-R\'{e}nyi
\frenchspacing
\newtheorem{theorem}{Theorem}[section]
\newtheorem{corollary}{Corollary}[section]
\newtheorem{conjecture}{Conjecture}[section]
\newtheorem{proposition}{Proposition}[section]
\newtheorem{lemma}{Lemma}[section]
\newtheorem{definition}{Definition}[section]
\newtheorem{example}{Example}[section]
\newtheorem{remark}{Remark}[section]
\newtheorem{case}{Case}[section]
\newtheorem{condition}{Condition}[section]
\newcommand{\pf}{\noindent {\it Proof:} }
\newcommand{\proof}{\noindent {\it Proof:} }

\title{{\bf\Large Bounded size biased couplings, log concave distributions and concentration of measure for occupancy models}}
\author[1]{Jay Bartroff}
\author[1]{Larry Goldstein}
\author[2]{\"{U}mit I\c{s}lak}
\affil[1]{University of Southern California} \affil[2]{Bo\u{g}azi\c{c}i University}

\footnotetext{This work was partially supported by NSA grant
H98230-11-1-0162 (Bartroff and Goldstein) and NSF grant DMS-1310127
(Bartroff).}

%Larry Goldstein and \"{U}mit
%I\c{s}lak\\University of Southern California} \
 \maketitle
\date{}

\begin{abstract}
Threshold-type counts based
on multivariate occupancy models with log concave marginals admit
bounded size biased couplings under weak conditions, leading to new
concentration of measure results for random graphs, germ-grain
models in stochastic geometry  and multinomial allocation models. The results obtained compare favorably with classical methods, including the use of McDiarmid's inequality,  negative association, and self bounding functions.
\end{abstract}

\section{Introduction}

A random graph on $m$ vertices in which edges are independently
 present between every two distinct vertices is one
framework that leads to an {\em occupancy model} described by a
 vector 
 \begin{equation}\label{occ.model}
{\bf M}=(M_\alpha)_{\alpha \in [m]}
\end{equation}
of
nonnegative integer valued random variables~$M_\alpha$, where $[m]=\{1,\ldots,m\}$ and $M_\alpha$ is the degree of
vertex~$\alpha$.  In such models, given a
nonnegative integer threshold $d \ge 0$, many authors have studied
the distribution of quantities such as \bea \label{intro:YgeY=d}
Y_{ge} = \sum_{\alpha \in [m]} {\bf 1}(M_\alpha \ge d) \qmq{and}
Y_{eq} = \sum_{\alpha \in [m]} {\bf 1}(M_\alpha =d) \ena which, in
the Erd\H{o}s-R\'{e}nyi random graph case just described, count the
number of vertices that have degree at least and exactly $d$,
respectively. Interest in the distributions of the random variables
defined in (\ref{intro:YgeY=d}) focuses on their approximation by
distributional limits such as the normal, and their finite sample
concentration properties. The purpose of the current manuscript is the latter, the study of the concentration of such random variables via the use of size biased couplings derived from the Stein's method literature. The
concentration of measure phenomenon has received a great deal of
attention since  the groundbreaking work of Talagrand \cite{Ta95}, and has found applications in areas as diverse as
statistics, random matrix theory, combinatorics, information theory,
and randomized algorithms. We refer to  \cite{Bo13} and \cite{Le05} for excellent treatments  of the subject. The results in this paper hold for more general occupancy models~\eqref{occ.model}, not just random graphs, and to refer to a generic occupancy model~\eqref{occ.model} we will call $M_\alpha$ a ``count'' and $\alpha$ an ``urn''.

To give the flavor of our results, continue to consider the
Erd\H{o}s-R\'{e}nyi random graph on $m$ vertices where each disjoint
pair of vertices is independently connected by an edge with
probability $p \in (0,1)$ and $M_\alpha$ is the degree of
vertex~$\alpha$.  The work \cite{GhGo11c} derived concentration
results for the number of isolated vertices, or equivalently,  for
the variable~$Y_{ge}$ in \eqref{intro:YgeY=d} with $d=1$. Here we allow each vertex to have its own threshold~$d_\alpha$ to either meet, exceed, or differ from, and which are allowed to take any value, each 
pair of disjoint edges
$\{i,j\}$ is to have its own connection probability
$p_{\{i,j\}}$, and each vertex to be
weighted according to a nonnegative `importance factor' 
$w_\alpha$.
In Theorem \ref{thm:ergraph}, for random graph models including the Erd\H{o}s-R\'{e}nyi, we provide sub-Poisson concentration bounds for random variables of the form \bea
\label{Y.ge.eq:def} Y_{ge} = \sum_{\alpha \in [m]} w_\alpha {\bf
1}(M_\alpha \ge d_\alpha) \qmq{and} Y_{ne} = \sum_{\alpha \in [m]}
w_\alpha {\bf 1}(M_\alpha \not = d_\alpha), \ena that is, for the
weighted number of components of ${\bf M}$ having size at least
$d_\alpha$, and not equal to $d_\alpha$, respectively. In addition, 
\begin{enumerate}
\item For the germ-grain models in stochastic geometry introduced in Section~\ref{sec:GG}, Theorems \ref{thm:GG} and \ref{thm:GGNN} provide sub-Poisson concentration results for the volume covered by multi-way intersections, and counts of neighbors, respectively.

\item For the multinomial model introduced in Section~\ref{mult:app}, Theorem \ref{thm:multocc} provides sub-Poisson concentration results for urn occupancy counts.

\end{enumerate}

The current work springs from that of \cite{GhGo11a} and \cite{GhGo11b},
which demonstrated how bounded size bias couplings can be used to
achieve concentration of measure results. Those works in turn were
built on the base of \cite{Ch07}, which showed how tools from
Stein's method (see \cite{St72} and \cite{St86}, and \cite{CGS10}
and \cite{Ro11} for overviews), and in particular exchangeable
pairs, can be used to expand the scope of application of the
concentration of measure phenomenon.  Through the use of bounded
size bias couplings, \cite{GhGo11b} produced concentration results
for examples including the number of relatively ordered subsequences
of a random permutation, the number of local maxima of a random
function on a lattice, the number of urns containing exactly one
ball in a uniform urn allocation model, and the volume covered by
the union of $n$ balls placed uniformly over a subset of
$\mathbb{R}^p$ with volume~$n$. In \cite{GhGo11c}, a concentration
result was obtained for the number of isolated vertices in the
Erd\H{o}s-R\'{e}nyi random graph.

Lemma~\ref{bd.szbc.d}, one main result in the present work, provides
a framework for the construction of bounded size bias couplings for
threshold counts of random variables having a discrete log concave
distribution. Such constructions allow the results of \cite{GhGo11b} and \cite{GhGo11c} to be extended to 
counts of multinomial urn occupancies that exceed or meet any
values, to the covered volume of multi-way intersections in germ-grain models, and to counts
of the number of vertices of the Erd\H{o}s-R\'{e}nyi  graph having any degrees.
Further, we do not require
an identical distribution assumption and consider occupancy
thresholds and importance weightings that may depend on the
component $\alpha \in [m]$. In Section~\ref{TalAzu} we show how our
results improve on what can be obtained by competing methods.

In order to prepare for the rest of the paper,  we provide some background  on size biased distributions and couplings. 
First,  recall that for a nonnegative random variable $Y$ with finite
positive mean $\mu$, we say that $Y^s$ has the $Y$-size bias
distribution if
\begin{equation}
\label{def:sb} E[Yf(Y)] =\mu E[f(Y^s)]
\end{equation}
for all functions $f$ for which these expectations
exist. For a survey on the diverse appearances of size biasing in
probability and statistics, see \cite{AGK}.  We say we have a size
bias coupling when a random variable $Y^s$ satisfying \eqref{def:sb}
is defined on the same space as $Y$, and the coupling is said to be
bounded when there exists $c \in [0,\infty)$ such that $|Y^s-Y| \le
c$ almost surely. The work \cite{GhGo11b} showed that for nonnegative~$Y$ with finite mean~$\mu$ and bounded size bias coupling~$Y^s$ satisfying $|Y^s-Y|
\le c$, if $Y^s\ge Y$ then 
\begin{equation}
 \label{a}
P\left(Y-\mu \le -t\right)\le \exp\left(-\frac{t^2}{2c \mu}\right)
\quad \mbox{for all $t > 0$.}
\end{equation}
And if the moment generating
function $m(\theta)=E(e^{\theta Y})$ is finite at $\theta=2/c$, then
\begin{equation}
\label{b} P\left(Y-\mu \ge t\right)\le
\exp\left(-\frac{t^2}{2c\mu+ct}\right) \qmq{for all $t>0$.}
\end{equation}

The bound  \eqref{a} holds without the monotonicity assumption and we prove this in Section~\ref{proof:main}, thus providing a left tail bound for any
application in \cite{GhGo11b} which previously lacked one. After a version of this manuscript was circulated, Theorem
\ref{thm:main} below of \cite{arba} removed the monotonicity
assumption using different methods, and further improved the result
of \cite{GhGo11b} by removing the assumption that the moment
generating function of $Y$ be finite at $2/c$, relaxing the bounded
coupling condition to $Y^s - Y \leq c$, and by improving the
inequality to \eqref{arba:bound} which, as shown there, implies
\eqref{a} and \eqref{b}. 

The subsequent work \cite{CGJ17} strictly generalizes Theorem 1.1 of \cite{arba} by, in Theorems 3.3 and 3.4, relaxing the almost sure boundedness assumption by the condition that there exists $p \in (0,1]$ such that $P[X^s \le X+c|X^s \ge x] \ge p$ for all $x$ for an upper tail bound, and $P[X^s \le X+c|X \le x] \ge p$ for all $x$ for the lower tail. Theorem 3.4 is in the spirit of Bennett's inequality, with upper bounds given in terms of a variance proxy, rather than the mean.

\begin{theorem}
\label{thm:main} Let $Y$ be a nonnegative random variable with
nonzero, finite mean $\mu$, and suppose there exists a coupling of
$Y$ to a variable $Y^s$ having the $Y$-size bias distribution that
satisfies $Y^s \le Y + c$ for some $c>0$ with probability one. Then
\begin{equation}
\label{arba:bound} \max\left\{\sup_{t \ge 0}P(Y-\mu \ge
t),\sup_{-\mu \le t \le 0}P(Y-\mu \le t) \right\} \le \left(
\frac{\mu}{\mu+t}\right)^{(t+\mu)/c}e^{t/c}.
\end{equation}
\end{theorem}
Note that the upper tail inequality given  in (\ref{arba:bound}) can
be rewritten in the more familiar form
$$
P(Y-\mu \geq t) \leq \exp \left(- \frac{\mu}{c}h\left(\frac{t}{\mu}
\right) \right) \qmq{for all $t>0$}
$$
where $h(x)= (1+x) \log(1+x)-x$, $x \geq -1$. Using the inequality
$$h(x) \geq \frac{x^2}{2+2x/3},\quad x\geq 0$$ (for example, see
\cite[Exercise 2.8]{Bo13}), one immediately obtains the following
Bernstein type inequality as a corollary, which provides a slight
improvement over (\ref{b}). 

\begin{corollary}% \label{bernstein}
In the setting of Theorem~\ref{thm:main}, 
\begin{equation}
\label{bernsteincorollary} P(Y- \mu \geq t) \leq \exp \left(-
\frac{t^2}{2c\mu +2ct/3} \right) \qmq{for all $t>0$.} 
\end{equation}
\end{corollary}

Next we briefly review constructions of random variables
having the size bias distribution of linear
combinations of indicator random variables.  Throughout we will write
${\cal L}(\cdot)$ for law, or  distribution.

We start by stating Lemma 4.1 of \cite{GoPe10}. When $A$ is an event satisfying $0<P(A)<1$ and ${\cal
F}$ is a $\sigma$-algebra, a simple application of nested
conditioning shows that for all bounded continuous functions $f$,
the random variable\footnote{In \eqref{sb:cond} and below, for an
event~$A$ we write ${\cal L}(Y|A)$ to denote the law of the random
variable with distribution $P(Y\in B|A)$. We abuse notation in the
standard way by writing  $P(\cdot|{\cal F})$ and $P(\cdot|R)$ to
denote conditioning on a $\sigma$-algebra ${\cal F}$ and the
$\sigma$-algebra $\sigma(R)$ generated by a random variable~$R$.} $Y=P(A|{\cal F})$
satisfies
\begin{equation}\label{sb:cond}
E(Yf(Y))=E(f(Y){\bf 1}_A) = \mathbb{P}(A) \mathbb{E}(f(Y )|A) \qmq{and hence} {\cal
L}(Y|A) ={\cal L}(Y^s).
\end{equation}

Next, Lemma~\ref{sblem} is a special case of a result
of \cite{GoRi96} that suggests constructions of size biased
couplings for sums of nonnegative random variables with finite
means.

\begin{lemma} \label{sblem}
Let $Y=\sum_{\alpha \in [m]}w_\alpha X_\alpha$ be a sum of Bernoulli
variables $(X_\alpha)_{\alpha \in [m]}$ weighted by nonnegative
constants $(w_\alpha)_{\alpha \in [m]}$ and satisfying $EY>0$.
Suppose that  for each $\alpha \in [m]$ the variables
$\{X_\beta^\alpha, \beta \in [m] \}$ are defined on a common probability space such
that \bea \label{Xbetaalphaiscond} {\cal L}(X_\beta^\alpha, \beta
\in [m])={\cal L}(X_\beta, \beta \in [m]|X_\alpha=1). \ena Then  for
each $\alpha \in [m]$ letting \bea \label{def:Yalpha}
Y^\alpha=\sum_{\beta \in [m]} w_\beta X_\beta^\alpha, \ena and $I$ a
random index with distribution \bea \label{def:Iprop}
P(I=\alpha)=\frac{w_\alpha EX_\alpha}{EY}, \ena the law ${\cal
L}(Y^I)$ given by the mixture $\sum_{\alpha \in [m]}
P(I=\alpha){\cal L}(Y^\alpha)$ is the $Y$-size bias distribution.
\end{lemma}

We note that $(X_\alpha)_{\alpha \in [m]}$ is allowed to have any
joint distribution with Bernoulli marginals. In addition, the joint
distributions of $\{X_\beta^\alpha, \beta \in [m]\}$ are constrained
only to satisfy \eqref{Xbetaalphaiscond}, and in particular, if all
variables are defined on a common space then the dependence
structure between the collections $(X_\beta^\alpha)_{\beta \in [m]}$
over $\alpha \in [m]$ can be arbitrary. In this latter case, and
when $I$ is defined on the common space, we construct $I$
independent of $(X_\beta^\alpha)_{\beta \in [m]}$ so that
$Eg(Y^I)=\sum_{\alpha \in [m]}P(I=\alpha)
E(g(Y^\alpha)|I=\alpha)=\sum_{\alpha \in
[m]}P(I=\alpha) E (g(Y^\alpha))$,  thus achieving the desired mixture of the lemma.

To understand the connections between \eqref{sb:cond} and
Lemma~\ref{sblem} let us briefly explain how the former implies the
latter. Suppose $Y$ is given as the weighted sum of Bernoulli
variables, as in Lemma~\ref{sblem}. Then letting
$w=\sum_{\alpha\in[m]} w_\alpha$, for an index  $J$ with
distribution $P(J=\alpha)=w_\alpha/w$, $\alpha \in [m]$, chosen
independently of $(X_\alpha)_{\alpha \in [m]}$, and $A=\{X_J=1\}$
and ${\cal F}=\sigma\{X_\alpha,\alpha \in [m]\}$ we obtain
\begin{equation}
\label{wP=Y} wP(A|{\cal F})=\sum_{\alpha \in [m]} w_\alpha X_\alpha
= Y.
\end{equation}
Taking expectation in (\ref{wP=Y}) yields $wP(A)=EY$. Now, if ${\cal
L}(Y')={\cal L}(Y|A)$ then, with $Y^\alpha$ as in (\ref{def:Yalpha})
and $I$ with distribution (\ref{def:Iprop}), the reader can easily
check that $E[g(Y^I)]=E[g(Y')]$ for all bounded continuous functions
$g$ so that, by \eqref{sb:cond}, ${\cal L}(Y^I)$ of
Lemma~\ref{sblem} is the $Y$-size bias distribution.

\begin{comment}
To give an idea of the flavor of our results, again consider the
Erd\H{o}s-R\'{e}nyi random graph on $m$ vertices where each disjoint
pair of vertices is independently connected by an edge with
probability $p \in (0,1)$. Let the component~$M_\alpha$ of the
vector ${\bf M}=(M_\alpha)_{\alpha \in [m]}$ record the degree of
vertex~$\alpha$. The work \cite{GhGo11c} derived concentration
results for the number of isolated vertices, or equivalently,  for
the variable ~$Y_{ge}$ in \eqref{intro:YgeY=d} with $d=1$. Here we
consider the random graph where each pair of disjoint edges
$\{i,j\}$ is allowed to have its own connection probability
$p_{\{i,j\}}$. Next, we allow each vertex $\alpha$ to have its own
threshold~$d_\alpha$ to either meet, exceed, or differ from, which
is allowed to take any value. Lastly, we allow each vertex to be
weighted according to a nonnegative `importance factor,' or weight,
$w_\alpha$. Thus, in the Erd\H{o}s-R\'{e}nyi random graph, and more
generally, we provide sub-Poisson concentration bounds of the form
\eqref{arba:bound} for random variables of the form \bea
\label{Y.ge.eq:def} Y_{ge} = \sum_{\alpha \in [m]} w_\alpha {\bf
1}(M_\alpha \ge d_\alpha) \qmq{and} Y_{ne} = \sum_{\alpha \in [m]}
w_\alpha {\bf 1}(M_\alpha \not = d_\alpha), \ena that is, for the
weighted number of components of ${\bf M}$ having size at least
$d_\alpha$, and not equal to $d_\alpha$, respectively.
\end{comment}

The rest of the paper is organized as follows. Section
\ref{sec:mainresults} shows how to construct bounded size bias
couplings for random variables of the form (\ref{Y.ge.eq:def}) when
the components $M_\alpha$ of ${\bf M}$ have a discrete log concave
distribution, with support bounded from below in the case of
 $Y_{ge}$, and when for all $\alpha \in [m]$ the remaining counts conditioned
on $M_\alpha=a$ and $M_{\alpha} = b$ for $a,b \in {\cal S}_\alpha$
can be closely coupled whenever $a$ is `close' to $b$.
In Section~\ref{applications} we provide complete descriptions of
the three models mentioned above, and apply the results of
Section~\ref{sec:mainresults} to obtain concentration of measure inequalities. A comparison of the size bias method  for
concentration with other techniques in the literature is included in
Section~\ref{TalAzu}. Section~\ref{proof:main} contains the proof that \eqref{a} holds without the monotonicity assumption $Y^s\ge Y$.

\section{Bounded coupling constructions under log concavity}\label{sec:mainresults}The purpose of this section is to form the theoretical background for size biased coupling
constructions that are to be used for obtaining concentration of measure inequalities for the statistics described in the Introduction.  First, we note that Lemma~\ref{sblem}  gives a recipe for the construction of a
variable having the $Y$-size biased distribution, and in particular
only suggests how a coupling may be created. Here, we construct couplings
not directly on the occupancy vectors ${\bf M}=(M_\alpha)_{\alpha
\in [m]}$ themselves, but on a collection of `more basic' variables ${\cal U}$ that we term {\em configurations}, and of which the occupancy counts are
functions. The configuration ${\cal U}$ will be specified for each application.
For instance, when ${\bf M}$ is the count of vertex
degrees in an Erd\H{o}s-R\'{e}nyi graph on a vertex set $[m]$ we
take the configuration ${\cal U}$ to be the collection of edge
indicator variables $(X_{\{\alpha,\beta\}})_{\alpha \not = \beta,
\{\alpha,\beta\} \subset [m]}$, and, similarly, when ${\bf M}$
counts the number of balls in each urn in a multinomial model, the
configuration ${\cal U}$ records the location of each ball.
In these examples the variables making up the
configurations~${\cal U}$ are random, but we also will refer to realizations of ${\cal U}$ as configurations, so configurations may also contain deterministic variables.

\begin{definition}\label{def:config}
When the occupancy counts ${\bf M}$ are given as $F({\cal U})$ for
some collection of (possibly random) variables ${\cal U}$ and measurable function $F$,
we say that ${\bf M}$ \textit{corresponds to} the configuration
${\cal U}$ through $F$,  or that configuration ${\cal U}$
\textit{has corresponding occupancy counts} ${\bf M}$ with respect
to $F$.
\end{definition}

    The function $F$ will be fixed in each of our applications, and so we may
     omit its mention when there is no possibility of confusion. In
     such cases we may write, for instance, that ${\cal U}$ has corresponding occupancy counts ${\bf M}$.

Specializing to the case of interest here, Lemma~\ref{sblem}
suggests the following size bias coupling construction for sums of the form
\eqref{Y.ge.eq:def}, say $Y_{ge}$ for concreteness. Following
\eqref{Xbetaalphaiscond}, given a configuration ${\cal U}$ from the
model, one constructs ${\cal U}^\alpha$, one for each $\alpha \in
[m]$, on the same space as ${\cal U}$, with law given by
\begin{equation}
\label{U.is.conditional.ge.dalpha} {\cal L}({\cal
U}^\alpha) = {\cal L}({\cal
U}|X_\alpha=1)={\cal L}({\cal U}|M_\alpha \ge d_\alpha).
\end{equation}
One
then obtains the variable $Y_{ge}^\alpha$ by evaluating the sum
$Y_{ge}$ of \eqref{Y.ge.eq:def} on the occupancy counts
corresponding to ${\cal U}^\alpha$, and the size bias variable
$Y_{ge}^s$ by selecting $Y^\alpha$ with probability proportional to
the expectation of $w_\alpha X_\alpha= w_\alpha {\bf 1}(M_\alpha \ge
d_\alpha)$, independently of all else.

For the construction of a configuration satisfying
\eqref{U.is.conditional.ge.dalpha}, Lemma~\ref{bd.szbc.d} and
Corollary~\ref{cor:all.at.once} below show how to achieve a bounded
coupling between $M_\alpha$ and a variable with distribution ${\cal
L}(M_\alpha|M_\alpha \ge d_\alpha)$ when the distribution of
$M_\alpha$ is log concave. Lemma~\ref{add.subtract.one} will be used
to construct the remainder of a configuration that has the correct
conditional counts for the urns $\beta \not = \alpha$ when their
marginals have the distribution of the sum of independent Bernoulli
variables, that is, when they have a Poisson Binomial distribution
\cite{Hoeffding56}; see \eqref{def:PB-dist} for a formal definition.
\label{p:PoisBin}

For any nonnegative integer $m$ let $[m]_0=\{0,\ldots,m\}$, and for
any subset ${\cal S}$ of $\mathbb{R}$ and any $t_1,t_2 \in
\mathbb{R}$, let $t_1{\cal S}+t_2=\{t_1s+t_2: s \in {\cal S}\}$. For
instance $[m]-1=[m-1]_0$ when $m \ge 1$. For a discrete random variable $M$ let
$p_x=P(M=x)$ and $\supp(M)=\{x\in\mathbb{R}: p_x>0\}$ be the
probability mass function and support of $M$, respectively. Recall
that $M$ is called a \textit{lattice} random variable if $\supp(M)
\subset h_1\mathbb{Z}+h_2$ for some real numbers $h_1\ne 0$, $h_2$.
We can without loss of generality assume our lattice random
variables $M$ have $\supp(M)\subset\mathbb{Z}$ by applying the
transformation
 $(M-h_2)/h_1$. Such a lattice random variable~$M$ is \textit{log concave~(LC)} if $\supp(M)$ is an \textit{integer interval}; that is, if
\beas %\label{def:int.int}
\supp(M)=(k_1,k_2)\cap \mathbb{Z}\qmq{for some} k_1,k_2\in
\mathbb{Z}\cup\{\pm\infty\},\; k_1<k_2-1, \enas  and
\begin{equation}\label{pLC}
p_x^2\ge p_{x-1}p_{x+1}\qmq{for all}x\in\mathbb{Z}.
\end{equation}

Under a lattice log concave assumption on the distribution of $M$,
Parts~\ref{pi.part} and \ref{rho.part} of Lemma~\ref{bd.szbc.d}
provide bounded couplings of random variables with distributions
${\cal L}(M|M\ge d)$ and ${\cal L}(M|M\le d)$, respectively, to
variables with distributions ${\cal L}(M|M\ge d+1)$ and ${\cal
L}(M|M\le d-1)$. Part~\ref{X.ne.d.part} shows that there is a
bounded coupling of $M$ to a variable having distribution ${\cal
L}(M|M \ne d)$, provided $M$ is not degenerate at $d$. These results
are extensions of \cite[Lemma 3.3]{GoPe10}, that showed the $d=1$
case of Part~\ref{pi.part} when $M$ is $\mbox{Bin}(n,p)$ with $p \in
(0,1)$.

In the following we let $\mbox{Bern}(p)$ denote the Bernoulli
distribution giving mass $1-p$ and $p$ to $0$ and $1$, respectively.

\begin{lemma} \label{bd.szbc.d}
 Let $M$ be a lattice LC random variable
with support ${\cal S} \subset \mathbb{Z}$.
\begin{enumerate}
\item\label{pi.part}   For $x, d \in \mathbb{Z}$ define
$$\pi_x^{(d)}=
\begin{cases}
\frac{P(M\ge x+1)P(M=d)}{P(M\ge d+1)P(M=x)},& \mq{if}x, d+1\in{\cal S}\qmq{and}x\ge d\\
0,&\mbox{otherwise.}
\end{cases}$$ Then the following hold.
\begin{enumerate}
\item\label{0<pi<1} $0\le\pi_x^{(d)}\le 1$ for all $x, d$.

\item\label{d->d+1} If $d+1\in{\cal S}$ and $N, Z$ are random variables such that ${\cal L}(N)={\cal L}(M|M\ge d)$ and ${\cal L}(Z|N)=\mbox{Bern}(\pi_N^{(d)})$, then ${\cal L}(N+Z)={\cal L}(M|M\ge d+1)$.
\end{enumerate}

\item\label{rho.part}   For $x, d \in \mathbb{Z}$ define
$$\rho_x^{(d)}=
\begin{cases}
\frac{P(M\le x-1)P(M=d)}{P(M\le d-1)P(M=x)},&  \mq{if}x, d-1\in{\cal S}\qmq{and}x\le d\\
0,&\mbox{otherwise.}
\end{cases}$$ Then the following hold.
\begin{enumerate}
\item\label{0<rho<1} $0\le\rho_x^{(d)}\le 1$ for all $x, d$.

\item If $d-1\in{\cal S}$ and $N, Z$ are random variables such that ${\cal L}(N)={\cal L}(M|M\le d)$ and ${\cal L}(Z|N)=\mbox{Bern}(\rho_N^{(d)})$, then ${\cal L}(N-Z)={\cal L}(M|M\le d-1)$.
\end{enumerate}

\item\label{X.ne.d.part}  Fix $d\in\mathbb{Z}$ such that $P(M=d)<1$. Let $Z_+, Z_-$ be conditionally independent given $M$ with ${\cal L}(Z_+|M)=\mbox{Bern}(\pi_M^{(d)})$ and ${\cal L}(Z_-|M)=\mbox{Bern}(\rho_M^{(d)})$.  Let $Z$ be independent of $Z_+$, $Z_-$, and $M$ with ${\cal L}(Z)=\mbox{Bern}(q)$, where
\begin{equation*}
q=\frac{P(M\ge d+1)}{P(M\ne d)}.
\end{equation*} Then
\begin{equation}\label{X.ne.d}
{\cal L}\left(M+X\right)={\cal L}(M|M\ne d),
\end{equation}
where  $X=ZZ_+-(1-Z)Z_-$.

\end{enumerate}
\end{lemma}

In other words, the conclusion \eqref{X.ne.d} says that, given $M$,
a random variable with distribution ${\cal L}(M|M\ne d)$ can be
formed by flipping an independent $q$-coin~$Z$  and, if heads,
adding $1$ to $M$ with probability $\pi_M^{(d)}$, and otherwise
subtracting $1$ with probability $\rho_M^{(d)}$. We note that when
$M<d$ (resp.\ $M>d$), the probability~$\pi_M^{(d)}$  of adding
(resp.\ $\rho_M^{(d)}$ of subtracting)  $1$ is $0$, and when $M=d$,
$M$ is changed with probability $1$ by either adding or subtracting
$1$. We also note that when $M$ achieves the upper or lower limit of
its support,  the Bernoulli probability of adding to, or subtracting
from $M$, respectively, is zero.

We define the hazard function of a lattice random variable $M$ with
support ${\cal S}$ as
\begin{equation}\label{haz}
h_x=\frac{P(M=x)}{P(M\ge x)}=\frac{p_x}{\sum_{y\ge x}p_y}\qmq{for}
x\in{\cal S}.
\end{equation}

To prove Lemma~\ref{bd.szbc.d} we require the following fact that lattice LC distributions have nondecreasing hazard functions. This is well known for continuous LC distributions, e.g., \cite{BaBe05,Muller02}.

\begin{lemma}\label{lem:haz}
If $M$ is lattice LC with support ${\cal S}$ then the hazard
function $h_x$ given in (\ref{haz}) is nondecreasing on ${\cal S}$.
\end{lemma}

\noindent {\em Proof:} For any $x,y\in {\cal S}$ with $x \le y$ note
that by (\ref{pLC}) we have
\begin{equation*}
\frac{p_{x+1}}{p_x}\ge \frac{p_{x+2}}{p_{x+1}}\ge\cdots\ge
\frac{p_{y+1}}{p_{y}}.
\end{equation*}
If $x, x+1\in {\cal S}$ then
\begin{multline*}
1/h_x-1/h_{x+1}=\sum_{y \in {\cal S}:\;y\ge x}p_y/p_x-\sum_{y \in {\cal S}:\;y\ge x+1}p_y/p_{x+1}=\sum_{y \in {\cal S}:\;y \ge x}(p_y/p_x-p_{y+1}/p_{x+1})\\
 =\sum_{y \in {\cal S}:\;y \ge x}\frac{p_y}{p_{x+1}}\left(\frac{p_{x+1}}{p_x}-\frac{p_{y+1}}{p_y}\right)\ge 0.
\end{multline*}\bbox

\bigskip

\noindent {\em Proof of Lemma~\ref{bd.szbc.d}:} Clearly $\pi_x^{(d)}
\ge 0$, and to show that $\pi_x^{(d)}\le 1$ it suffices to assume
that $d, d+1\in {\cal S}$ since $\pi_x^{(d)}=0$ otherwise. Let $h_x$
be the hazard function of $M$ defined by \eqref{haz}. For any $d\le
x\in {\cal S}$, by Lemma~\ref{lem:haz} we have $h_d \le h_x$, and
therefore
$$\pi_x^{(d)}=\frac{1/h_x-1}{1/h_d-1}\le 1,$$
proving Part~\ref{0<pi<1}.

To prove Part~\ref{d->d+1}, letting $p_x=P(M=x)$ and $G_x=P(M\ge
x)$, for any $k=1,2,\ldots$ we have
\begin{multline*}
P(N+Z\ge d+k)=P(N\ge d+k)+P(N=d+k-1, Z=1)\\
=P(M\ge d+k|M\ge d)+\pi_{d+k-1}^{(d)}P(M=d+k-1|M\ge d)
=\frac{G_{d+k}}{G_{d}}+\left(\frac{G_{d+k}p_d}{G_{d+1}p_{d+k-1}}\right)\frac{p_{d+k-1}}{G_{d}} \\
 =\frac{G_{d+k}}{G_{d}G_{d+1}}\left(G_{d+1}+p_d
\right) =\frac{G_{d+k}}{G_{d}G_{d+1}} G_{d}=\frac{G_{d+k}}{G_{d+1}}
=P(M\ge d+k|M\ge d+1).
 \end{multline*}

For Part~\ref{0<rho<1} let $\wtilde{M}=-M$, which is LC.  For
$d-1\in {\cal S}$ and $d\ge x\in {\cal S}$,
$$
\rho_x^{(d)}=\frac{P(\wtilde{M}\ge
-x+1)P(\wtilde{M}=-d)}{P(\wtilde{M}\ge
-d+1)P(\wtilde{M}=-x)}=\wtilde{\pi}_{-x}^{(-d)} \in [0,1]
$$
by Part~\ref{0<pi<1}, where $\wtilde{\pi}$ is defined with respect
to $\wtilde{M}$. The rest of the proof of Part~\ref{rho.part} is
similar to that of Part~\ref{pi.part}.

Moving to Part~\ref{X.ne.d.part}, letting $N$ denote the random
variable $M+X$ on the LHS of \eqref{X.ne.d}, we will show that
\begin{equation}\label{Ltail.N0} P(N\le y)=P(M\le y|M\ne d)\qmq{for
all}y\in\mathbb{Z},\; y<d,
\end{equation}
the proof that $P(N\ge y)=P(M\ge y|M\ne d)$ for all $y>d$ being
similar. Fix $y<d$ and without loss of generality assume that
\begin{equation}\label{d-k+1inS}
y+1\in {\cal S},
\end{equation}
 since otherwise \eqref{Ltail.N0} holds trivially as both sides are $0$ or $1$. With $p_x, G_x$ as above and
$F_x=P(M\le x)$,
 \begin{align}
P(N\le y)&=P(M\le y-1)+P(M=y,\; Z=0)+P(M=y,\; Z=1,\; Z_+=0)\nonumber\\
&+P(M=y+1,\; Z=0,\; Z_-=1)\nonumber\\
&=F_{y-1}+p_{y}(1-q)+p_{y}q(1-\pi_{y}^{(d)})+p_{y+1}(1-q)\rho_{y+1}^{(d)}\nonumber\\
&=F_{y}+ p_{y+1}(1-q)\rho_{y+1}^{(d)},\label{Ltail.N}
\end{align} this last because $\pi_{y}^{(d)}=0$ since $y<d$. If $d-1\in{\cal S}$ then \eqref{Ltail.N} is
$$F_{y}+ p_{y+1}\left(1-\frac{G_{d+1}}{1-p_d}\right) \left( \frac{F_{y}p_d}{F_{d-1}p_{y+1}} \right) =F_{y}+\left(\frac{F_{d-1}}{1-p_d}\right)
\frac{F_{y}p_d}{F_{d-1}} =\frac{F_{y}}{1-p_d} =P(M\le y| M\ne d).$$
Otherwise $d-1\not\in{\cal S}$ so $\rho_{y+1}^{(d)}=0$, hence
\eqref{Ltail.N} is $F_{y}$.   If $y=d-1$ then $\min{\cal S}=d$ by
virtue of the assumption \eqref{d-k+1inS}, so $$P(N\le
d-1)=F_{d-1}=0=P(M\le d-1|M\ne d).$$ In the remaining case,
$d-1\not\in{\cal S}$ and $y\le  d-2$, we have $\max{\cal S}<d-1$
again by virtue of \eqref{d-k+1inS}, and in particular
$d\not\in{\cal S}$. Then
$$P(M\le y|M\ne d)=P(M\le y)=F_{y}=P(N\le y),$$ finishing the proof.\bbox

\begin{corollary} \label{cor:all.at.once}
Let $M$ be a lattice LC random variable with support~${\cal S}$
satisfying $a:=\inf {\cal S}>-\infty$, and let $d \in {\cal S}$.
Then one can construct a random variable $A$ on the same space as
$M$ such that ${\cal L}(M+A)={\cal L}(M|M \ge d)$ and $0 \le A \le
d-a$.
\end{corollary}

\noindent {\em Proof:} It suffices to prove the $a=0$ case  because
given $M$ satisfying the hypotheses and $d\in{\cal S}$ we have
$d':=d-a\in\mbox{supp}(M-a)$. Using that $M-a$ is an LC random
variable, the $a=0$ case of the corollary guarantees $0\le A\le d'$
such that
$${\cal L}(M-a+A)={\cal L}(M-a|M-a \ge d') = {\cal L}(M-a|M \ge d).$$ Using this fact, we have
$$P(M+A\le x)=P(M-a+A\le x-a)=P(M-a\le x-a|M\ge d)=P(M\le x|M\ge d),$$ which is the desired result.

To prove the $a=0$ case we successively construct random variables
$M_0,\ldots,M_d$, all on the same space as $M$, such that
\begin{equation}
\label{eq:inductMk} {\cal L}(M_k)={\cal L}(M|M \ge k) \qmq{for
$k=0,1,\ldots,d$.}
\end{equation}
Letting $M_0=M$, \eqref{eq:inductMk} is satisfied for $k=0$. For
$k=0,\ldots,d-1$, given $M_0,\ldots,M_k$ satisfying the
distributional equality in  \eqref{eq:inductMk}, let $X_k$ be a
Bernoulli random variable on the same space as $M_0,\ldots,M_k$
satisfying \beas {\cal L}(X_k|M_k) = \mbox{Bern}(\pi_{M_k}^{(k)})
\qmq{and set} M_{k+1}=M_k+X_k. \enas
%Clearly $M_0$ satisfies
%\eqref{eq:inductMk} for $k=0$. Assume that $M_k$ satisfies
%\eqref{eq:inductMk} for some $k=0,\ldots,d-1$.
It is easily checked that ${\cal L}(M)$ being lattice LC implies
${\cal L}(M_k)$ is LC. Hence Part~\ref{d->d+1} of
Lemma~\ref{bd.szbc.d} yields that $M_{k+1}$ satisfies
\eqref{eq:inductMk} for $k+1$. In particular, for $k=d$ our
construction yields $M_d=M+A$ with $A=X_0+\ldots+X_{d-1}$ satisfying
$0 \le A \le d$, and \eqref{eq:inductMk} yields the desired
distributional property, concluding the proof. \bbox

\bigskip

 Corollary~\ref{cor:all.at.once}, a consequence of
    Part~\ref{pi.part} of  Lemma~\ref{bd.szbc.d}, shows the
    existence of a `uniformly close coupling' of an LC `urn count' random variable $M$ to one with
    distribution ${\cal L}(M|M\ge d)$, and will be applied to coupling constructions for $Y_{ge}$. Similarly, Part~\ref{X.ne.d.part} of Lemma~\ref{bd.szbc.d}, depending on Parts~\ref{pi.part} and \ref{rho.part}, will be applied to $Y_{ne}$.

Recall that a random variable $M$  is said to have a Poisson
Binomial distribution with parameter ${\bf p}=(p_j)_{j \in [m]}$,
denoted by $M \sim {\cal PB}({\bf p})$, when
 \begin{equation}
\label{def:PB-dist} {\cal L}(M)= {\cal L}\left(\sum_{j \in [m]}
 X_j\right)
\end{equation}
where $X_j$ are independent Bernoulli random variables with
$P(X_j=1)=p_j$ for $j \in [m]$.

When there exists $p$ such that $p_j=p$ for all $j \in [m]$,
then $ M \sim \mbox{Bin}(m,p)$. We note that the
distribution of a single Bernoulli random variable, with support
$\{0,1\}$, trivially satisfies (\ref{pLC}) and hence is LC. Since
\cite{Keilson71} demonstrates that LC is preserved under
convolution, the claim of the following lemma is immediate.
\begin{lemma} \label{PB.is.LC}
The Poisson Binomial distribution ${\cal PB}({\bf p})$ is LC.
\end{lemma}

%\begin{comment}
%Turning to the hypergeometric distribution, it is shown in
%\cite[Theorem~A]{Ehm91} that a hypergeometric random variable can be
%written as a sum of independent but non-identically-distributed
%Bernoulli random variables, hence the following lemma is a special
%case of the previous.
%\begin{lemma} \label{Hy.is.LC}
%The hypergeometric distribution is LC.
%\end{lemma}
%\end{comment}

When $M$ has distribution ${\cal PB}({\bf p})$ for ${\bf p}=(p_j)_{j
\in [m]}$, then for all $d \in \mathbb{Z}$ we have \beas
P(M=d)=q_{eq}(d,{\bf p}) \qmq{where}  q_{eq}(d,{\bf p})=\sum_{s
\subset [m],\; |s|=d} \prod_{j \in s}p_j \prod_{j \not \in s}(1-p_j)
\enas and so \beas P(M \ge d)=q_{ge}(d,{\bf p}) \qmq{and} P(M \not =
d)=q_{ne}(d,{\bf p}) \enas where
\begin{equation}
\label{def:qge.qeq}
q_{ge}(d,{\bf p})=\sum_{k=d}^m q_{eq}(k,{\bf p}) \qmq{and}
q_{ne}(d,{\bf p}) = 1- q_{eq}(d,{\bf p}).
\end{equation}

With $\alpha \in [m]$, the majority of our constructions
make use of the following definition for collections of counts of the form
\begin{equation}
\label{N.props}
\{\mathbf{N}_a^\alpha = (N_{\beta,a}^\alpha)_{\beta\in[m]}:
a\in\mS_\alpha\},
\end{equation}
where $\mS_\alpha$, $\alpha \in [m]$, are given support sets. Below, for given weight ${\bf w}=(w_\alpha)_{\alpha \in [m]}$ and threshold
${\bf d}=(d_\alpha)_{\alpha \in [m]}$ vectors, let 
\begin{equation}
\label{def.absw.absd} |{\bf w}|=\max_{\alpha \in [m]} w_\alpha
\qmq{and} |{\bf d}|=\max_{\alpha \in [m]} d_\alpha.
\end{equation}

\begin{definition}\label{def:Bstara} For $B \ge 0$  and $\alpha \in [m]$, we say that the
collection of counts \eqref{N.props}  has Property $(B,\ge,\alpha)$ if
\begin{equation}\label{PropAgeB}
\sum_{\beta \not = \alpha} w_{\beta} {\bf
        1}(N_{\beta,a+1}^\alpha \ge d_\beta) \le \sum_{\beta \not = \alpha}
    w_{\beta} {\bf 1}(N_{\beta,a}^\alpha \ge d_\beta)+|{\bf w}|B \qm{for all
        $\{a,a+1\} \subset {\cal S}_\alpha$,}
\end{equation}
and Property $(B,\not =, \alpha)$ if
\begin{equation}
\label{PropABnot} \sum_{\beta \not = \alpha} w_{\beta} {\bf
        1}(N_{\beta,b}^\alpha \not = d_\beta) \le \sum_{\beta \not = \alpha}
    w_{\beta} {\bf 1}(N_{\beta,a}^\alpha \not = d_\beta)+|{\bf w}|B \qm{for all
        $\{a,b\} \subset {\cal S}_\alpha$ with $|b-a|=1$,}
\end{equation}
for all ${\bf w} =(w_1,\ldots,w_m) \in (0,\infty)^m$ and all ${\bf
d}=(d_1,\ldots,d_m) \in \mathbb{Z}^m$. If the counts  have Property $(B, \ge, \alpha)$  (resp.  $(B,\not =, \alpha)$) for all $\alpha \in [m]$, then they are said to have Property $(B, \ge)$ (resp. $(B, \not =)$). 

For $\star \in \{ \ge , \not = \}$ we
say that a collection of configurations has Property $(B,\star)$
when their corresponding occupancy counts do.
\end{definition}

The following claims are immediate.
\begin{enumerate}

\item For given $\alpha\in[m]$, counts \eqref{N.props}
  have Property $(B,\ge,\alpha)$ and $(B,\not =,\alpha)$  if
\bea \label{bound.difference.by.B} \Bvert \{\beta: \beta \not =
\alpha, N_{\beta,a}^\alpha \not = N_{\beta,a+1}^\alpha\}\Bvert \le B
\qm{ whenever $\{a,a+1\} \in {\cal S}_\alpha$.} \ena

\item For given $\alpha\in[m]$, counts \eqref{N.props} 
have Property $(0,\ge,\alpha)$  when
\bea \label{B.NA.case} N_{\beta,a+1}^\alpha \leq
N_{\beta,a}^\alpha\qm{ for all $\{a,a+1\} \subset {\cal S}_\alpha$
and $\beta \neq \alpha$.} \ena 

\item For given $\alpha\in[m]$, if counts \eqref{N.props} 
have Property $(B,\ge,\alpha)$ then
\bea \label{atobcriterion}\sum_{\beta \not = \alpha}w_{\beta}{\bf
1}(N_{\beta,b}^\alpha \ge d_\beta) \le \sum_{\beta \not = \alpha}
w_{\beta} {\bf 1}(N_{\beta,a}^\alpha \ge d_\beta)+ |{\bf w}|B(b-a)
\ena for all $\{a,b\} \subset {\cal S}_\alpha$ with $a\leq b$, all
${\bf w} \subset (0,\infty)^m$ and ${\bf d} \in \mathbb{Z}^m$.
\end{enumerate}

Now let ${\cal U}$ be a configuration corresponding to an
occupancy model ${\bf M}=(M_\alpha)_{\alpha \in [m]}$, and for all
$\alpha \in [m]$ let ${\cal S}_\alpha$ be the support of $M_\alpha$,
$a_\alpha=\inf {\cal S}_\alpha$, and $b_\alpha = \sup {\cal
S}_\alpha$. For $a \in {\cal S}_\alpha$ let
\begin{equation}
\label{ind.add.one} {\cal L}({\cal V}_a^\alpha) := {\cal L}({\cal
U}|M_\alpha =a).
\end{equation}

Theorem~\ref{+1-1:generally} is our main tool for the construction of bounded size biased couplings
	for $Y_{ge}$ and $Y_{ne}$ for all applications other than those in Section~\ref{sec:GG}. All constants $B$ in the following, whose value may change between different occurrences, are universal.

\begin{theorem} \label{+1-1:generally}
Let ${\cal U}$ be a configuration corresponding to occupancy counts
${\bf M}=(M_\alpha)_{\alpha \in [m]}$, where for all $\alpha \in
[m]$ the component $M_\alpha$ is lattice LC.  Suppose that  for all  $\alpha \in [m]$ there exists configurations 
$\{{\cal U}_a^\alpha, a \in {\cal S}_\alpha\}$ on a common
space satisfying
\begin{equation}
\label{calU=dcalVa} {\cal L}({\cal U}_a^\alpha)= {\cal L}({\cal
V}_a^\alpha) \qm{for all $a \in {\cal S}_\alpha$,}
\end{equation}
where ${\cal V}_a^\alpha$ is given by~\eqref{ind.add.one}.
\begin{enumerate}
\item If $\{{\cal U}_a^\alpha, a\in\mS_\alpha\}$ has Property $(B,\ge)$ and $a_\alpha>-\infty$ for all $\alpha \in [m]$, then there exists a coupling of variables $Y$ and $Y^s$ on the same space such that ${\cal L}(Y)={\cal L}(Y_{ge})$ and ${\cal L}(Y^s)={\cal L}(Y_{ge}^s)$
satisfying
\begin{equation}\label{lem-bound-Yge}
Y^s \le Y + |{\bf w}|\left(B |{\bf d}-{\bf a}| + 1\right),
\end{equation}
where ${\bf d}-{\bf a}=(d_\alpha-a_\alpha)_{\alpha \in [m]}$.

\item If $\{{\cal U}_a^\alpha, a\in\mS_\alpha\}$ has Property $(B,\ne)$ for all $\alpha \in [m]$, then there exists a coupling of variables $Y$ and $Y^s$ on the same space such that ${\cal L}(Y)={\cal L}(Y_{ne})$ and ${\cal L}(Y^s)={\cal L}(Y_{ne}^s)$ satisfying
\begin{equation}
\label{lem-bound-Yne} Y^s \le Y + |{\bf w}|\left(B + 1\right).
\end{equation}
\end{enumerate}
\end{theorem}

\noindent {\em Proof of Theorem \ref{+1-1:generally}:} We prove \eqref{lem-bound-Yge} first. Fix
$\alpha \in [m]$ and let $N_\alpha$ with distribution ${\cal
L}(M_\alpha)$ be defined on the same space as, and independent of,
the configurations $\{{\cal U}_a^\gamma, \gamma \in [m], a \in {\cal
S}_\gamma \}$. By Corollary~\ref{cor:all.at.once} one can construct
$A_\alpha$ on the same space as $N_\alpha$ such that ${\cal
L}(N_\alpha+A_\alpha)={\cal L}(M_\alpha|M_\alpha \ge d_\alpha)$ with
$0 \le A_\alpha \le d_\alpha-a_\alpha$. In particular,
\begin{multline} \label{UNalphaXalpha.Uconditioned}
{\cal L}({\cal U}|M_\alpha \ge d_\alpha) = \sum_{a \ge d_\alpha, a \in {\cal S}_\alpha}{\cal L}({\cal V}_a^\alpha)P(M_\alpha=a|M_\alpha \ge d_\alpha)\\
%= \sum_{a \ge d_\alpha}{\cal L}({\cal %V}_a^\alpha)P(N_\alpha=a|N_\alpha \ge d_\alpha)
= \sum_{a \ge d_\alpha,a \in {\cal S}_\alpha}{\cal L}({\cal
V}_a^\alpha)P(N_\alpha+A_\alpha=a)={\cal L}({\cal
V}_{N_\alpha+A_\alpha}^\alpha)={\cal L}({\cal
U}_{N_\alpha+A_\alpha}^\alpha).
\end{multline}
Clearly ${\cal L}({\cal U}_{N_\alpha}^\alpha)={\cal L}({\cal U})$ by
 \eqref{ind.add.one} and
\eqref{calU=dcalVa}. Let ${\bf N}^\alpha$ and ${\bf N}_{ge}^\alpha$
be the counts corresponding to ${\cal U}_{N_\alpha}^\alpha$ and
${\cal U}_{N_\alpha+A_\alpha}^\alpha$, respectively. Since, for all
$\alpha \in [m]$, the configurations $\{{\cal U}_a^\alpha, a \in
{\cal S}_\alpha \}$ have Property $(B,\ge,\alpha)$ and $0 \le A_\alpha \le
d_\alpha-a_\alpha$, we have
\begin{multline}
\label{+1bd.alpha} \sum_{\beta \in [m]} w_{\beta} {\bf
1}(N_{\beta,ge}^\alpha \ge d_\beta) - \sum_{\beta \in [m]}w_{\beta}
{\bf
1} (N_{\beta}^\alpha \ge d_\beta) \le |{\bf w}| (B  A_{\alpha} + 1)\\
\leq |{\bf w}| ( B(d_{\alpha} - a_{\alpha}) + 1) \le |{\bf w}| (
B|{\bf d} - {\bf a}| + 1),
\end{multline}
where we applied observation \eqref{atobcriterion}, and where the
factor $+1$ accounts for the maximum possible change from $0$ to $1$
of the indicator associated to urn $\alpha$.

Let $I$ be a random index with distribution \eqref{def:Iprop}
defined with respect to the weighted indicators in
\eqref{Y.ge.eq:def} summing to $Y_{ge}$, independent of all other
variables, and set ${\bf N}={\bf N}^I$ and ${\bf N}_{ge}={\bf
N}_{ge}^I$. Now, for nonnegative integer counts ${\bf
n}=(n_\alpha)_{\alpha \in [m]}$, setting \bea \label{Ybfn.ge} Y({\bf
    n})=\sum_{\alpha \in [m]} w_\alpha {\bf 1}(n_\alpha \ge d_\alpha),
\ena averaging \eqref{+1bd.alpha} over $\alpha$ distributed as $I$
we obtain \beas  Y({\bf N}_{ge}) \le Y({\bf N})+|{\bf w}| (B|{\bf
d}-{\bf a}|+1). \enas

The counts ${\bf N}$ have distribution ${\cal L}({\bf M})$, as the
same holds for ${\bf N}^\alpha$ for all $\alpha \in [m]$, by virtue
of ${\cal L}({\cal U}_{N_\alpha}^\alpha)={\cal L}({\cal U})$. In
particular $Y=Y({\bf N})$ has distribution ${\cal L}(Y_{ge})$. By
\eqref{UNalphaXalpha.Uconditioned} the indicators
$X_\beta^\alpha={\bf 1}(N_{\beta,ge}^\alpha \ge d_\beta)$ satisfy
\eqref{Xbetaalphaiscond} with $X_\beta={\bf 1}(M_\beta \ge
d_\beta)$, and Lemma~\ref{sblem} yields $Y^s=Y({\bf N}_{ge})$ has
the $Y_{ge}$-size biased distribution. The proof of
\eqref{lem-bound-Yge} is now complete.

To prove \eqref{lem-bound-Yne}, first recall that we have reduced to
the case that  $P(M_\alpha \not = d_\alpha)<1$ for all $\alpha \in
[m]$, allowing us to invoke Part~\ref{X.ne.d.part} of
Lemma~\ref{bd.szbc.d}. Construct $Z_-,Z_+$ and $Z$ as in the lemma,
on the same space as $M_\alpha$, so that with
$A_\alpha=ZZ_+-(1-Z)Z_-$ we have ${\cal L}(M_\alpha+A_\alpha)={\cal
L}(M_\alpha|M_\alpha \not = d_\alpha)$, with $-1 \le A_\alpha \le
1$. The proof proceeds in the same  way as for
\eqref{lem-bound-Yge}.

Let ${\bf N}^\alpha$ and ${\bf N}_{ne}^\alpha$ be the counts
corresponding to ${\cal U}_{N_\alpha}^\alpha$ and ${\cal
U}_{N_\alpha+A_\alpha}^\alpha$, respectively. Since, for all $\alpha
\in [m]$, the configurations $\{{\cal U}_a^\alpha, a \in {\cal
S}_\alpha \}$ have Property $(B,\ne,\alpha)$ and $-1 \le A_\alpha \le 1$,
we have \bea \label{ne.bd.alpha} \Bvert \sum_{\beta \in [m]}
w_{\beta} \left( {\bf 1}(N_{\beta,ne}^\alpha \not = d_\beta) -  {\bf
1} (N_{\beta}^\alpha \not = d_\beta) \right)\Bvert \le |{\bf w}|
(B|A_{\alpha}|+1)  \le |{\bf w}| (B+1). \ena

Let $I$ be an independent index with distribution \eqref{def:Iprop}
defined with respect to the weighted indicators in
\eqref{Y.ge.eq:def} summing to $Y_{ne}$ and set ${\bf N}={\bf N}^I$
and ${\bf N}_{ne}={\bf N}_{ne}^I$. Then with $Y({\bf n})$ given by
\eqref{Ybfn.ge} with $\ge$ replaced by $\not =$, we have ${\cal
L}(Y({\bf N}))={\cal L}(Y_{ne})$ and, by Lemma~\ref{sblem}, that
${\cal L}(Y({\bf N}_{ne}))={\cal L}(Y_{ge}^s)$. Now averaging
\eqref{ne.bd.alpha} over $I$ yields \beas |Y({\bf N}_{ne})-Y({\bf
N})| \le |{\bf w}|(B+1), \enas and the desired conclusion.
\bbox
\bigskip

The following lemma is helpful in verifying that the conditions of
Theorem~\ref{+1-1:generally} are in force when the configurations
${\cal U}$ corresponding to the occupancy counts ${\bf M}$ are given
in terms of independent Bernoulli variables.

\begin{lemma} \label{add.subtract.one}
    Let ${\cal X}=(X_\alpha)_{\alpha \in [m]}$ be a collection of independent Bernoulli random variables with
    respective success probabilities $p_1,\ldots,p_m\in (0,1)$, and let
    $R=\sum_{\alpha \in [m]} X_\alpha$.
%For $a \in [m]$ let
%\beas
%{\cal L}(Q_a) = {\cal L}(R||R|=a).
%\enas
Then there exists $\{{\cal X}_a, a\in[m]_0 \}$ defined on a common
space such that, for $a \in [m-1]_0$, \bea \label{U.coupled.in.a}
{\cal L}({\cal X}_a)={\cal L}({\cal X}|R=a) \qmq{and} {\cal X}_a \le
{\cal X}_{a+1} \qm{with probability one,} \ena where for $\{{\bf
x},{\bf y}\} \subset \{0,1\}^m$ we write ${\bf x} \le {\bf y}$ when
$x_i \le y_i$ for all $i \in [m]$.
\end{lemma}

\noindent {\em Proof:} Recall that the density $p(n)$ of an integer
valued random variable is a P\'olya frequency function of order 2
(or simply, is ${\rm PF}_2$; see \cite{Schoenberg51}) when\beas
\Bvert
\begin{array}{cc}
p(m_1-n_1) & p(m_1-n_2) \\
p(m_2-n_1) & p(m_2-n_2)
\end{array}
\Bvert \ge 0 \qmq{for all $m_2 \ge m_1$ and $n_2 \ge n_1$.} \enas
For $\phi(y_1,\ldots,y_m)$ a coordinatewise non-decreasing function
of $(y_1,\ldots,y_m) \in \mathbb{R}^n$, \cite{Efron65} Section 3,
shows that if $Y_1,\ldots,Y_m$ are independent integer valued random
variables with ${\rm PF}_2$ densities then
$E(\phi(Y_1,\ldots,Y_m)|Y_1+ \cdots+Y_m=a)$ is a non-decreasing
function of $a$. As the Bernoulli density is ${\rm PF}_2$, we find
in particular that \bea \label{a.monotone}
E(\phi(X_1,\ldots,X_m)|R=a) \le E(\phi(X_1,\ldots,X_m)|R=a+1)
\qm{for all $a \in [m-1]_0$}, \ena for $\phi$ a coordinatewise
non-decreasing function on $\{0,1\}^m$. See \cite{BrPeRo} for a
simple proof of this fact in the Bernoulli case. 

Relation \eqref{a.monotone} is expressed in Definition 2.1 of
\cite{liggett} as ${\cal X}_a \le {\cal X}_{a+1}$ with probability
one. Hence by Theorem 2.4 of \cite{liggett} there exists a
distribution $Q_a(\cdot,\cdot)$  on $\{0,1\}^m \times \{0,1\}^m$ for
$({\cal V},{\cal W})$ such that \bea \label{for.prop.2} {\cal
L}({\cal V})={\cal L}({\cal X}_a), \quad {\cal L}({\cal W})={\cal
L}({\cal X}_{a+1}) \qmq{and} {\cal V} \le {\cal W} \qm{with
probability one.} \ena (See \cite{BrPeRo} for a specific
construction of the pair $({\cal V},{\cal W})$.)

With some slight abuse of notation, let $Q_a(\cdot|\cdot)$ denote
the $Q_a(\cdot,\cdot)$ conditional distribution of the second
argument given the first. Let ${\cal X}_0$
be the vector in $\{0,1\}^m$ with all coordinates equal to 0, and
for $a=1,\ldots,m$, given ${\cal X}_{a-1}$ let ${\cal X}_a$ be
sampled from the conditional distribution $Q_a(\cdot|{\cal
X}_{a-1})$. Clearly ${\cal L}({\cal X}_a)={\cal L}({\cal X}|R=a)$
for $a=0$. Assuming this identity holds for $a\in [m-1]_0$, it holds
also for $a+1$, as ${\cal L}({\cal X}_{a+1})$ is the conditional law
$Q_{a+1}(\cdot|{\cal X}_a)$ averaged over the distribution ${\cal
L}({\cal X}|R=a)$, which equals ${\cal L}({\cal X}|R=a+1)$ by
construction.

Hence the first property in \eqref{U.coupled.in.a} holds; the second
is a consequence of the last relation in \eqref{for.prop.2}. \bbox

\section{Applications} \label{applications}
We now present in detail the three models mentioned in the Introduction, and use the
constructions in Section~\ref{sec:mainresults} to prove
concentration bounds for each case. With the exception of the volume
of multi-way intersections in germ-grain models,  the variables of
interest are weighted occupancy counts of the form
\begin{equation}\label{Y.occ}
Y_{ge} = \sum_{\alpha \in [m]} w_\alpha {\bf 1}(M_\alpha \ge
d_\alpha) \qmq{and} Y_{ne} = \sum_{\alpha \in [m]} w_\alpha {\bf
1}(M_\alpha \not = d_\alpha),
\end{equation}
although see the next paragraph for related random variables that can also be handled. Without loss of generality we assume that all summands in
\eqref{Y.occ} are non-constant, for if a summand were constant then it could simply be subtracted from the corresponding $Y$ and the number of summands~$m$ decremented by one. Some consequences of this assumption are that all $w_\alpha$ are strictly
positive and, with ${\cal S}_\alpha$ denoting
the support of $M_\alpha$, that $\inf
{\cal S}_\alpha < d_\alpha < \sup {\cal S}_\alpha+1$  for all $\alpha \in [m]$ when
considering $Y_{ge}$, and $0<P(M_\alpha \not = d_\alpha)<1$  for all $\alpha \in [m]$ when
considering $Y_{ne}$. 

The concentration bounds we provide for variables of the form
(\ref{Y.occ}) also yield bounds for the `complementary' sums
\beas %\label{beta.minus}
\sum_{\alpha \in [m]} w_\alpha {\bf 1}(M_\alpha < d_\alpha) =
\sum_{\alpha \in [m]}w_\alpha  - Y_{ge} \qmq{and}
 \sum_{\alpha \in [m]} w_\alpha
{\bf 1}(M_\alpha = d_\alpha) =\sum_{\alpha \in [m]}w_\alpha  -
Y_{ne} , \enas with the mean~$\mu=EY$ replaced by $\sum_{\alpha \in
[m]} w_\alpha -\mu$ and the roles of the right and left tails
reversed. In fact, all our results can be extended further, with
essentially only a notational burden, to random variables of the
form \beas Y = \sum_{\alpha \in [m]} w_\alpha {\bf 1}(M_\alpha
\star_\alpha d_\alpha), \qmq{where} \star_\alpha \in \{\ge, \not
=\}, \enas and therefore, in like manner, to the sums of
complementary form.

Lastly, we note that when $M_\alpha \sim {\cal PB}({\bf p}_\alpha)$ for each $\alpha \in
[m]$, by
    \eqref{def:qge.qeq} the means $\mu_{ge}$ and $\mu_{ne}$ of
$Y_{ge}$ and $Y_{ne}$ are given respectively by
\begin{equation}
\label{means.multinomial.occupancy} \mu_{ge}=\sum_{\alpha \in [m]}
w_\alpha q_{ge}(d_{\alpha},{\bf p}_\alpha) \qmq{and}
\mu_{ne}=\sum_{\alpha \in [m]} w_\alpha q_{ne}(d_{\alpha},{\bf
p}_\alpha),
\end{equation}
where $q_{ge}$ and $q_{ne}$ are given by \eqref{def:qge.qeq}.

\subsection{Degree counts in Erd\H{o}s-R\'{e}nyi type graphs}\label{sec:E-R}
The classical Erd\H{o}s-R\'{e}nyi random graph on $m$ vertices is
constructed by placing an edge between each pair of distinct
vertices independently and with equal probability. The model was
originally used in conjunction with the probabilistic method for
proving the existence of graphs with certain properties (see
\cite{alon}) and has been popular more recently for modeling complex
networks (e.g.,  \cite{chung}).

The classical
Erd\H{o}s-R\'{e}nyi graph with constant connectivity~$p$ has
been the object of much study. Asymptotic normality of the number of
vertices of degree $d$ was shown in \cite{karonski} when
$m^{(d+1)/d} \rightarrow \infty$ and $mp \rightarrow 0,$ or $m p
\rightarrow 0$ and $mp -\log m -d\log \log m \rightarrow -\infty.$
Asymptotic normality when $mp \rightarrow c >0$ was obtained in
\cite{Barbour89}. Optimal bounds in the Kolmogorov metric can be
found in \cite{Korkecki90} and \cite{Go12}. Other univariate results
on asymptotic normality of counts on random graphs are given in
\cite{janson}, and references therein. Smooth function bounds were
obtained in \cite{GoRi96} for the vector whose $k$ components count
the number of vertices of fixed degrees $d_1, d_2, \ldots, d_k$ when
$p = \theta/(m - 1) \in (0,1)$ for fixed $\theta$, implying
asymptotic multivariate joint normality. This work was later
extended in \cite{LiRe12} to the inhomogeneous random graph model
which will be the setting in the current paper.

 Here we consider graph degree counts when the likelihood of an edge may depend on the identity of the vertices it connects. Formally, let ${\cal G}_m$ be an Erd\H{o}s-R\'{e}nyi random graph on
the vertices $[m]$, where the presence of an edge joining  distinct
vertices $\alpha$ and $\beta$ is recorded by the indicator
    $X_{\{\alpha,\beta\}}$ with success probability
$p_{\{\alpha,\beta\}}$, with all such indicators independent. We set
$p_{\{\alpha,\alpha\}}=0$ for all $\alpha \in [m]$, making all
$X_{\alpha,\alpha}$ identically zero. The classical model is
recovered by setting $p_{\{\alpha,\beta\}}=p$ for some $p \in [0,1]$
for all $\alpha \not = \beta$.

For $Y_{ge}$, with similar remarks applying to $Y_{ne}$, by removing
any edge $\{\alpha,\beta\}$ with $p_{\{\alpha,\beta\}}=1$ and
decrementing each of the two thresholds $d_{\alpha}, d_{\beta}$ by
one we may assume that $p_{\{\alpha,\beta\}}<1$ for all
$(\alpha,\beta) \in [m] \times [m]$. Having also reduced to the case
where all the indicators in (\ref{Y.occ}) are nontrivial allows us
to assume that $\sum_{\beta:\; \beta \neq \alpha}
p_{\{\alpha,\beta\}} > 0$ for all $\alpha \in [m].$

Let the components $M_\alpha$ of ${\bf M}=(M_\alpha)_{\alpha \in
[m]}$ record the degree of vertex $\alpha$, that is,
$M_\alpha=\sum_{\beta \in [m]}X_{\{\alpha,\beta\}}$. By the definition~\eqref{def:PB-dist}, $M_\alpha \sim {\cal PB}({\bf
p}_\alpha)$ with ${\bf p}_\alpha = (p_{\{\alpha,\beta\}})_{\beta:
\beta \not = \alpha}$. By \eqref{def:qge.qeq} the means $\mu_{ge}$
and $\mu_{ne}$ of $Y_{ge}$ and $Y_{ne}$ have the  form
\eqref{means.multinomial.occupancy}. With $|{\bf w}|$ and $|{\bf
d}|$ as in \eqref{def.absw.absd}, let
\begin{equation}\label{C.occ}
c_{ge}=|{\bf w}|(|{\bf d}|+1) \qmq{and} c_{ne}=2|{\bf w}|.
\end{equation}

\begin{theorem} \label{thm:ergraph}
Concentration of measure inequalities
(\ref{a})-(\ref{bernsteincorollary}) hold for
counts $Y_{ge}$ and $Y_{ne}$ given by \eqref{Y.occ} in ${\cal G}_m$ for all $m \ge 1$, with corresponding $\mu$ and $c$ given by \eqref{means.multinomial.occupancy} and \eqref{C.occ}.
\end{theorem}

Theorem \ref{thm:ergraph} is a direct consequence of the following lemma.

\begin{lemma}\label{lem:ergraph}
In ${\cal G}_m$ there exists a coupling of $Y_{ge}$ to $Y_{ge}^s$, having the
$Y_{ge}$-size biased distribution, that satisfies $Y_{ge}^s-Y_{ge}
\le c_{ge}$, and a coupling of $Y_{ne}$ to $Y_{ne}^s$, having the
$Y_{ne}$-size biased distribution, satisfying $Y_{ne}^s-Y_{ne} \le
c_{ne}$.
\end{lemma}

\noindent {\em Proof:} In this model we take the configuration \beas
{\cal U}=\{X_{\{\gamma,\delta\}}, \{\gamma,\delta\} \subset [m]\},
\enas 
the collection of the independent Bernoulli edge
indicator variables of the graph ${\cal G}_m$. As the corresponding
counts ${\bf M}$ have LC marginal distributions with supports ${\cal
S}_\alpha$ satisfying $\inf {\cal S}_\alpha=0$, in order to invoke
Theorem~\ref{+1-1:generally} it is only required to show that for all
$\alpha \in [m]$, configurations $\{{\cal U}_a^\alpha,a \in [m]\}$
exist with Properties $(1,\geq)$ and $(1,\neq)$ satisfying
\eqref{calU=dcalVa}.

With $a \in {\cal S}_\alpha$ and ${\cal V}_a^\alpha$ as in
\eqref{ind.add.one}, by independence we obtain
\begin{multline*}
{\cal L}({\cal V}_a^\alpha) = {\cal L}({\cal U}|M_\alpha=a)
={\cal L}(X_{\{\gamma,\delta\}},  \{\gamma,\delta\} \subset [m]|M_\alpha=a)\\
={\cal L}\left(X_{\{\gamma,\delta\}}, \{\gamma,\delta\} \subset [m]\Bvert\sum_{\delta \in [m]} X_{\{\alpha,\delta\}}=a\right) \\
= {\cal L}\left(X_{\{\alpha,\delta\}}, \delta \in [m] \Bvert
\sum_{\delta \in [m]} X_{\{\alpha,\delta\}}=a \right)\times {\cal
L}\left(X_{\{\gamma,\delta\}}, \{\gamma,\delta\} \not \ni \alpha
\right),
%\\
%:={\cal L}(X_{\{\alpha,\delta\},a}^\alpha, \delta \in [m]) %\times {\cal L}\left(X_{\{\gamma,\delta\}}, \alpha \not \in %\{\gamma,\delta\}\right),
\end{multline*}
where here $\times$ denotes product measure.

On the same space and independently of ${\cal U}$ and of each other
over $\alpha \in [m]$, let \beas {\cal X}_a^\alpha =
\{X_{\{\alpha,\delta\},a}^\alpha, \delta \in [m]\}, \qmq{$a \in
[m]$} \enas be the collections of Bernoulli variables guaranteed by
Lemma~\ref{add.subtract.one} when taking $(X_\alpha)_{\alpha \in
[m]}$ in the lemma to be $(X_{\{\alpha,\delta\}})_{\delta \in [m]}$.
In particular \beas {\cal L}(X_{\{\alpha,\delta\},a}^\alpha, \delta
\in [m])={\cal L}\left(X_{\{\alpha,\delta\}}, \delta \in [m] \Bvert
\sum_{\delta \in [m]} X_{\{\alpha,\delta\}}=a \right), \enas and
setting ${\cal U}_a^\alpha = \{X_{\{\alpha,\delta\},a}^\alpha,
\delta \in [m] ,X_{\{\gamma,\delta'\}},  \{\gamma,\delta'\} \not \ni
\alpha\}$,
 we obtain ${\cal L}({\cal U}_a^\alpha)={\cal L}({\cal V}_a^\alpha)$.

Again by Lemma~\ref{add.subtract.one}, with $a+1 \in {\cal
S}_\alpha$, the Bernoulli variables equal to one in the collection
${\cal X}_{a+1}^\alpha$ are those equal to one in ${\cal
X}_a^\alpha$, with one additional variable. The configurations
${\cal U}_a^\alpha$ and ${\cal U}_{a+1}^\alpha$ therefore correspond
to graphs ${\cal G}_{m,a}^{\alpha}$ and ${\cal G}_{m,a+1}^{\alpha}$,
where the edge set of the latter is that of the former plus exactly
one additional edge attached to vertex $\alpha$. In particular, the
corresponding counts ${\bf N}_a^{\alpha}$ and ${\bf
N}_{a+1}^{\alpha}$ agree in all but coordinate $\alpha$ and one
additional coordinate. Hence \eqref{bound.difference.by.B} is
satisfied with $B=1$, implying that the configurations $\{{\cal
U}_a^\alpha,a \in [m]\}$ have both Properties $(1,\geq)$ and
$(1,\neq)$, thus completing the proof.

\bbox

\bigskip

In the standard case of equal thresholds $d_\alpha =d$ and unit weightings,
the expectations \eqref{means.multinomial.occupancy} of $Y_{ge}$ and
$Y_{ne}$ simplify  to \bea \label{erstandardmeans} \mu_{ge} = m
P(\mbox{Bin}(m-1, p) \geq d) \qmq{and} \mu_{ne} =m P(\mbox{Bin}(m-1,
p) \neq d) , \ena respectively, and the bounds
\eqref{a}-\eqref{bernsteincorollary} apply to $Y_{ge}$ with $c =
d+1$, and for $Y_{ne}$ with $c=2$. In particular, (\ref{a}) and
(\ref{bernsteincorollary}) yield that, for all $t>0$,
\begin{multline} \label{sbERnonis}
P(Y_{ge}-\mu_{ge} \leq -t) \leq
\exp\left(-\frac{t^2}{2(d+1)\mu_{ge}}\right) \,\, \mbox{and} \\
P(Y_{ge}-\mu_{ge} \geq t) \leq
\exp\left(-\frac{t^2}{2(d+1)(\mu_{ge}+t/3)}\right).
\end{multline}

The special case of the number of isolated vertices \beas Y_{is} =
\sum_{\alpha \in [m]} \mathbf{1}(M_{\alpha}=0) \enas for the
standard Erd\H{o}s-R\'{e}nyi model was handled in \cite{GhGo11c},
using an unbounded size bias coupling, and with much greater effort.
Techniques of the present paper can be used to obtain concentration
bounds for $Y_{is}$ in a much simpler way by noting that
$m-Y_{is}=Y_{ge}$ under unit weightings and equal thresholds
$d_\alpha=1$. In particular, the bounds \eqref{sbERnonis} hold with
$Y_{is}-\mu_{is}$ replacing $Y_{ge}-\mu_{ge}$ and setting $d=1$,
reversing the roles of the left and right tail bounds, and replacing
$\mu_{ge}$ by $m-\mu_{is}$. The left tail bound obtained in this
fashion is stronger than the corresponding bound \beas
P(Y_{is}-\mu_{is} \leq -t) \leq \exp \left( -\frac{t^2}{4\mu_{is}}
\right), \enas given in \cite{GhGo11c}, for $t \leq 6m(1-p)^{m-1}
-3m$, with similar remarks applying to the right tail.

Although the unbounded size bias coupling argument given in
\cite{GhGo11c} applies only to the case of isolated vertices,
Theorem~\ref{thm:ergraph} applies equally
for all degrees $d$. In particular, keeping $p$ and $d$ fixed and
letting $m \rightarrow \infty$, the left and right tail bounds for
$Y_{ge}$ provided by (\ref{a}) and (\ref{bernsteincorollary}), say,
will behave as $\exp(-t^2/(2(d+1)m))$ and
$\exp(-t^2/(2(d+1)(m+t/3)))$, respectively.

A well studied asymptotic is the case where $d$ is fixed and $m p \rightarrow \lambda$ for some $\lambda > 0$ so that for large $m$, the distribution $Bin(m-1,p)$ is close to a Poisson random variable with
	parameter $\lambda$. Here focusing
on the statistic $Y_{ge}
= \sum_{\alpha \in [m]} \mathbf{1}(M_{\alpha} \geq 1)$ for
simplicity, the mean satisfies $\mu_{ge} \rightarrow m(1 -
e^{-\lambda})$, and the resulting left and right tail bounds are
asymptotic to $\exp(-t^2/(4 m(1- e^{-\lambda})))$ and $\exp(-t^2/(4
m(1- e^{-\lambda})+4t/3))$, respectively, as $m \rightarrow \infty$.
Comparisons of these tail bounds with other techniques from the
literature will be discussed in detail in Section~\ref{TalAzu}.

\subsection{Germ-Grain models}\label{sec:GG}
Germ-grain models consist of sets (grains) placed at centers (germs)
determined by a random point process in some multidimensional space.
These models are used in applications including forestry
\cite{sher}, material science \cite{ande} and wireless sensor
networks \cite{dous}. For concreteness, here we consider models on
 the space $C_n=[0,n^{1/p})^p \subset \mathbb{R}^p$, equipped with the
Euclidean toroidal distance~$D$.

In the models considered in this section, a configuration
\label{config.comment} ${\cal V}$ is given by a collection $(v_\alpha)_{\alpha \in [m]}$ of points in $C_n$, and each
point~$v_\alpha$ is associated with a closed ball~$B_\alpha$
centered at $v_\alpha$ with positive radius $\rho_\alpha$.  We let
${\cal U}$ consist of points $U_1,\ldots,U_m$ sampled independently
in $C_n$ with strictly positive densities $f_1(x),\ldots,f_m(x)$ in
$C_n$, respectively. The positivity condition on the densities are
assumed for convenience, as along with \eqref{n.large} in Subsection
\ref{subsec:vol}, and \eqref{wide.enough.box} in Subsection
\ref{subsec:number}, respectively, it implies that the support is
$[m]_0$ for the number~$M(x,{\cal U})$ of intersections of ${\cal
U}$ at a point~$x \in C_n$, defined below in \eqref{def:Mx}, and
support $[m-1]_0$ for the number~$M_\alpha$ of neighbors of $U_\alpha$,
$\alpha \in [m]$, defined below in \eqref{Ma.GG1}.

\subsubsection{Volume of multi-way intersections in germ-grain models}
\label{subsec:vol}  For a point $x \in C_n$ and a configuration
${\cal V},$ let 
\begin{equation}
\label{def:Mx} M(x,{\cal V})=\sum_{\alpha \in
[m]} {\bf 1}(x \in B_\alpha),
\end{equation}
the number of balls $B_\alpha$
that contain the point $x$ in the configuration ${\cal V}$. In this
subsection we make the assumption that $n$ is large enough that \bea
\label{n.large} \sqrt{p} n^{1/p} > 2 \sum_{\alpha \in
[m]}\rho_\alpha. \ena Under \eqref{n.large} there exist
$(v_\alpha)_{\alpha \in [m]} \subset C_n$ such that $B_\alpha \cap
B_\beta = \emptyset$ for all $\alpha \not = \beta$, implying the
support of $M(x,{\cal U})$ is $[m]_0$.

For a fixed measurable function $d:C_n \rightarrow [m]$ and $f \in
[m]$ let \bea \label{def:Bge} {\bf
    1}_{x,ge}({\cal V}) = {\bf 1}\left( x \in B_{ge}(d(x),{\cal V}) \right) \qmq{where}
B_{ge}(f,{\cal V})= \bigcup_{\stackrel{r \subset [m]}{|r| \ge f}}
\bigcap_{\alpha \in r}B_\alpha. \ena Hence, dropping the dependence
on ${\cal V}$ unless clarity demands it, ${\bf 1}_{x,ge}=1$ if and
only if $x$ is contained in at least $d(x)$ of the balls
$B_\alpha,\alpha \in [m]$. Further, emphasizing dependence on the
function $d$ by writing ${\bf 1}_{x,d}={\bf 1}_{x,ge}$, the
indicators \beas {\bf 1}_{x,eq} = {\bf 1}_{x,d}-{\bf 1}_{x,d+1}
\qmq{and} {\bf 1}_{x,ne} = 1-{\bf 1}_{x,eq} \enas take the value 1
if $x$ is, and is not, contained in exactly $d(x)$ of the balls
$B_\alpha, \alpha \in [m]$, respectively. In particular, given a
nonnegative, bounded function $w(x)$ over $C_n$,
\begin{equation}
\label{Y.GG2} Y_{ge}({\cal V}) = \int_{C_n} w(x) {\bf
1}_{x,ge}({\cal V})\, dx \qmq{and} Y_{ne}({\cal V}) = \int_{C_n}
w(x){\bf 1}_{x,ne}({\cal V})\, dx
\end{equation}
are the volumes, weighted by $w(x)$, of the collection of points $x$
in $C_n$ contained in at least $d(x)$ balls, and some number of
balls other than $d(x)$ respectively. When $w(x)=1$ the variables in
\eqref{Y.GG2} are the volumes of the sets of points $x \in C_n$ that
are part of $d(x)$ way intersections, and intersections of size
other than $d(x)$, respectively.

Letting  \bea \label{def:palpha} p_\alpha(x)=P(D(x,U_\alpha) \le
\rho_\alpha), \quad \alpha \in [m], \ena  the variable $M(x,{\cal
U})$ of \eqref{def:Mx} has the Poisson Binomial distribution ${\cal
PB}({\bf p}(x))$. As, for instance, $\{x \in B_{ge}(x,{\cal V})\} =
\{M(x,{\cal V}) \ge d(x)\}$, we may write \bea \label{Y.GG3}
Y_{ge}({\cal V}) = \int_{C_n} w(x) {\bf 1}(M(x,{\cal V}) \ge d(x))\,
dx \,\,\mbox{and}\,\,Y_{ne}({\cal V})= \int_{C_n} w(x){\bf
1}(M(x,{\cal V}) \not = d(x))\, dx \ena whose expectations when
${\cal V}={\cal U}$ are given respectively by \bea \label{mu.GG2}
\mu_{ge} = \int_{C_n} w(x) q_{ge}(d(x),{\bf p}(x))dx\,\,
\qmq{and}\,\, \mu_{ne} = \int_{C_n} w(x) q_{ne}(d(x),{\bf p}(x))dx.
\ena

The size biased couplings for the germ-grain models in this
subsection and the following one are simple extensions of those in
Section 4 of \cite{GoPe10}. There are three differences between the
present case and that of \cite{GoPe10}. First, in that previous work
one considers only $d(x)=1$ for all $x \in C_n$ for the covered
volume. Next, $U_1,\ldots,U_m$ were taken in \cite{GoPe10} to have
the uniform distribution over $C_n$, and lastly the radii
$\rho_\alpha$ were set to some fixed $\rho$ for all $\alpha \in
[m]$. The imposition of these conditions result in various
simplifications in \cite{GoPe10}, which the following outline
generalizes in these aspects.

Let $U_0$ be sampled independently of  ${\cal U}$ with density
$f_0(x)=w(x){\bf 1}(x \in C_n)/w$ where $w=\int_{C_n}w(x)dx$, and
for $u \in C_n$ consider the event $F(u) = \left\{ M(u,{\cal U}) \ge
d(u) \right\}$ that the point $u$ lies in a $d(u)$ way intersection.
Then, in view of \eqref{Y.GG3}, it is easy to see that $Y_{ge}({\cal
U})=wP(F(U_0)|{\cal F})$ where ${\cal F}$ is the $\sigma$-algebra
generated by $\{U_\alpha, \alpha \in [m]\}$. Hence, by
\eqref{sb:cond}, $E(Y_{ge}({\cal U})| F(U_0))$ has the $Y_{ge}({\cal
U})$ size biased distribution. In particular, letting $Y_{ge}^s$
denote a variable with the $Y_{ge}({\cal U})$-size biased
distribution, we have that \bea \label{evaluate.at.M0} {\cal
L}(Y_{ge}({\cal U}^0))={\cal L}(Y_{ge}^s) \qmq{where} {\cal L}({\cal
U}^0) = {\cal L}({\cal U}|F(U_0)), \ena that is, the law of ${\cal
U}^0$ is that of the configuration ${\cal U}$ conditioned on the
event that the additional randomly chosen point $U_0$ lies in a
$d(U_0)$ way intersection.

In the following, for $\alpha \in [m], u \in C_n$ and $i \in
\{0,1\}$ let $f_{\alpha,u,i}$ be the distributions specified by \bea
\label{f01.conditionals} f_{\alpha,u,0}(E)=P(U_\alpha \in E|
D(U_\alpha,u) > \rho_\alpha) \qmq{and} f_{\alpha,u,1}(E)=P(U_\alpha
\in E| D(U_\alpha,u) \le \rho_\alpha). \ena

With the given weight, threshold functions $w(x), d(x)$ and radii
$(\rho_\alpha)_{\alpha \in [m]}$, let \beas |{\bf w}|=\sup_{x \in
    C_n}|w(x)|\quad |{\bf d}|=\sup_{x \in C_n}|d(x)| \qmq{and}
|\bm{\rho}|=\max_{\alpha \in [m]}\rho_\alpha, \enas and with $\pi_p$
the volume of the unit ball in $\mathbb{R}^p$, let
\begin{equation}\label{C.GG2}
c_{ge}=\pi_p |{\bf w}||{\bf d}||\bm{\rho}|^p  \qmq{and} c_{ne}=
\pi_p |{\bf w}||\bm{\rho}|^p.
\end{equation}

\begin{theorem} \label{thm:GG}
Concentration of measure inequalities
(\ref{a})-(\ref{bernsteincorollary}) hold for all $m \ge 1$ for
the volume covered by multi-way intersections in the germ-grain model described above, with $Y$, $\mu$ and $c$ given by
\eqref{Y.GG2}, \eqref{mu.GG2} and \eqref{C.GG2}.
\setcounter{monsterThmCounter}{\theenumi} 
\end{theorem}

Again, Theorem \ref{thm:GG} follows from the following bounded coupling construction.

\begin{lemma} \label{VolinRp}
In the germ-grain model described above there exists a coupling of $Y_{ge}$ to $Y_{ge}^s$, having the
    $Y_{ge}$-size biased distribution, that satisfies $Y_{ge}^s \le
    Y_{ge} + c_{ge}$, and a coupling of $Y_{ne}$ to $Y_{ne}^s$, having
    the $Y_{ne}$-size biased distribution, satisfying $Y_{ne}^s \le
    Y_{ne} + c_{ne}$.
\end{lemma}

\noindent {\bf Proof:}   As for any event $E$ we have
    \beas
    P(E|F(U_0)) = \frac{P(E,F(U_0))}{P(F(U_0))}=\int_{C_n} \frac{P(E,F(u))}{P(F(u))} \frac{P(F(u))}{P(F(U_0))}\frac{w(u)}{w}du,
    \enas
    we see ${\cal L}({\cal U}|F(U_0))$ of \eqref{evaluate.at.M0} is the mixture
    \bea \label{size.bias.mixture}
    {\cal L}({\cal U}|F(U_0)) = \int_{C_n} {\cal L}({\cal U}|F(u)) {\widetilde w}(u)du \qmq{where} {\widetilde w}(u)=\frac{P(F(u))}{P(F(U_0))}\frac{w(u)}{w}.
    \ena
    Hence, we sample ${\widetilde U}$ on $C_n$ with density ${\widetilde w}(u)$, and for ${\widetilde U}=u$, construct a configuration ${\cal U}$ with law ${\cal L}({\cal U}|F(u))$ in order to achieve ${\cal L}({\cal U}|F(U_0))$;
    see Lemma 2.3 of \cite{AGK} regarding size biasing mixture distributions.

For $u \in C_n$ let $(X_{\alpha,u})_{\alpha \in [m]}$ be independent
Bernoulli variables with success probabilities
$(p_\alpha(u))_{\alpha \in [m]}$ as given in \eqref{def:palpha}, and
for $a \in [m]_0$, let
\begin{equation}
\label{def:Xalphaxa}
    P(X_{1,u,a}=i_1,\ldots,X_{m,u,a}=i_m) = P\left(X_{1,u}=i_1,\ldots,X_{m,u}=i_m\left|\sum_{\alpha \in [m]}i_\alpha=a\right.\right).
\end{equation}
and
$$
{\cal X}_a=\{X_{\alpha,u,a}, \alpha \in [m]\}, \qmq{$a \in [m]_0$}
$$
be the collection of Bernoulli variables guaranteed by
Lemma~\ref{add.subtract.one} for the vector of success probabilities
$(p_\alpha(u))_{\alpha \in [m]}$.

For $\alpha \in [m]$ let $U_{\alpha;0}$ and $U_{\alpha;1}$ have
distributions $f_{\alpha,u,0}$ and $f_{\alpha,u,1}$ respectively, as
specified in \eqref{f01.conditionals}. Then for $a \in [m]_0$,
writing $U_{\alpha,a}=U_{\alpha;X_{\alpha,u,a}}$ and letting
$B_{\alpha,a}$ be the ball of radius $\rho_\alpha$ centered at
$U_{\alpha,a}$, the configuration ${\cal U}_a = \{U_{\alpha,a},
\alpha \in [m] \}$ satisfies ${\cal L}({\cal U}_a) = {\cal L}({\cal
U}|M(u,{\cal U})=a)$.

Let $N_{u}$ be constructed on the same space with distribution
${\cal L}(M(u,{\cal U}))$. Applying Corollary~\ref{cor:all.at.once}
to $M(u,{\cal U})$ we obtain $A_{u}$ satisfying \beas {\cal
L}(N_{u}+A_{u}) = {\cal L}(M(u,{\cal U})|F(u)) \qmq{with} 0 \le
A_{u} \le d(u). \enas Denoting the underlying configurations
corresponding to $N_{u}$  and to $N_{u}+A_{u}$, by ${\cal
U}_{N_{u}}$ and ${\cal U}_{N_{u}+A_{u}}$, respectively, we have
${\cal L}({\cal U}_{N_{u}+A_{u}})={\cal
    L}({\cal U}|F(u))$. Since ${\tilde U}=u$ was sampled according to density ${\widetilde w}(u)$, by \eqref{size.bias.mixture} for the second identity, letting $N=N_{\tilde U}$ and $A=A_{\tilde U}$,
\beas {\cal L}({\cal U}_N)={\cal L}({\cal U}) \qmq{and} {\cal
L}({\cal U}_{N+A})={\cal L}({\cal U}|F(U_0)). \enas By
\eqref{evaluate.at.M0}, $Y_{ge}^s=Y_{ge}({\cal U}_{N+A})$ has the
$Y_{ge}=Y({\cal U}_N)$ size biased distribution.

Note that \eqref{U.coupled.in.a} of Lemma~\ref{add.subtract.one}
implies that for $a \in [m-1]_0$ the indicators ${\cal X}_{a+1}$ and
${\cal X}_a$ are equal, but for one index, say $\beta_a$, such that
$X_{\beta_a,U_0,a}=0$ and $X_{\beta_a,U_0,a+1}=1$. Hence the
configurations ${\cal U}_a$ and ${\cal U}_{a+1}$ are the same but
for the one point indexed by $\beta_a$. Therefore, with $B_{ge,a}$
given by \eqref{def:Bge} with $B_\alpha$ replaced by $B_{\alpha,a}$,
 \bea \label{geq.indicators} {\bf 1}_{x,ge}({\cal U}_{a+1})  \le
{\bf 1}_{x,ge}({\cal U}_a){\bf 1}(x \not \in B_{\beta_a}^{a+1})+{\bf
1}(x \in B_{\beta_a}^{a+1}) \le {\bf 1}_{x,ge}({\cal U}_a)  + {\bf
1}(x \in B_{\beta_a}^{a+1}), \ena implying by \eqref{Y.GG2} that
\bea \label{eq:b-a} Y_{ge}({\cal U}_b) = \int_{C_n} w(x) {\bf
1}_{x,ge}({\cal U}_b) dx \le Y_{ge}({\cal U}_a) + \pi_p |{\bf w}|
(b-a) |\bm{\rho}|^p \ena for $b=a+1$, and hence for all $0 \le a \le
b \le m$. In particular, \beas Y_{ge}^s =Y_{ge}({\cal U}_{N+A}) \le
Y_{ge}({\cal U}_N)+ \pi_p |{\bf w}| A |\bm{\rho}|^p \le Y_{ge}({\cal
U}_N)+  \pi_p |{\bf w}||{\bf d}||\bm{\rho}|^p = Y_{ge}+ \pi_p |{\bf
w}||{\bf d}||\bm{\rho}|^p, \enas thus verifying the claim for
$Y_{ge}$.

To handle $Y_{ne}$, with $N$ as before satisfying ${\cal L}(N)={\cal
L}(M(U_0,{\cal U}))$, construct $A$ as in Part~\ref{X.ne.d.part} of
Lemma~\ref{bd.szbc.d}, so that $-1 \le A \le 1$ and \beas {\cal
L}(N+A)={\cal L}(M(U_0)|M(U_0,{\cal U}) \not = d(U_0)). \enas
Arguing as before, we have that $Y_{ne}^s=Y_{ne}({\cal U}_{N+A})$
has the $Y_{ne}=Y({\cal U}_N)$-size biased distribution.  As the
configurations ${\cal U}_a$ and ${\cal U}_{a+1}$ differ only at the
point indexed by $\beta_a$, with $\Delta$ denoting the symmetric
difference of sets, as in \eqref{geq.indicators} one has \beas {\bf
1}_{x,ne}({\cal U}_a) \le {\bf 1}_{x,ne}({\cal U}_b)+{\bf 1}(x \in
B_{\beta_{a \wedge b}}^{a} \Delta B_{\beta_{a\wedge b}}^{b}) \qm{for
$|a-b| \le 1, \{a,b\} \subset [m].$}
\enas As $|A| \le 1$, \beas Y_{ne}^s =Y_{ne}({\cal U}_{N+A}) \le
Y_{ne}({\cal U}_N)+  \pi_p |{\bf w}| |A| |\bm{\rho}|^p  \le Y_{ne}+
\pi_p |{\bf w}||\bm{\rho}|^p, \enas thus verifying the claim for
$Y_{ne}$. \bbox

\subsubsection{Number of neighbors in germ-grain models}
\label{subsubneighborsggm}  \label{subsec:number} In this section we
consider the occupancy model ${\bf M}({\cal V})=(M_\alpha({\cal
V}))_{\alpha \in [m]}$, based on the configuration ${\cal
V}=(v_\alpha)_{\alpha \in [m]}$, with components
\begin{equation}
\label{Ma.GG1} M_\alpha ({\cal V})=  \sum_{\beta \in [m] \setminus
\{\alpha\}} {\bf 1}(B_\alpha \cap B_\beta \not = \emptyset),
\end{equation}
where $B_\alpha$ is the closed unit ball centered at $v_\alpha,
\alpha \in [m]$. That is, $M_\alpha$ counts the number of
neighbors of $v_\alpha$, where we say $v_\alpha$ and $v_\beta$ are neighbors, and write $v_\alpha \sim v_\beta$, when $\alpha \not =
\beta$ and $B_\alpha \cap B_\beta \not = \emptyset$. The variables $Y_{ge}$ and $Y_{ne}$ are
again given as in \eqref{Y.occ} for ${\cal V}={\cal U}$, and are the
weighted sums of the contributions from points $U_\alpha$ of ${\cal
U}$ that have at least $d_\alpha$ neighbors, and a number other than
$d_\alpha$ neighbors, respectively. We drop the dependence on ${\cal V}$ unless clarity demands it.

Suppressing the dimension $p$ in our notation, in this section we
specialize to the unit radius $\rho_\alpha=1$ case in order to allow
Lemma \ref{lem:kappa} below to yield a bound on our coupling in
simple terms of classical geometric constants related to the
`kissing numbers' $\kappa_1^*$; see \cite{CoSl99} and \cite{Zong99}.
With $B_0$ the closed unit ball of radius one centered at the
origin, the constant $\kappa_1^*$ is the maximum number of closed
unit balls in $\mathbb{R}^p$ that can be packed so that their
closures intersect $B_0$, with all balls having disjoint interiors.
The related constant $\kappa_1$ arises below, which is the maximum
number of unit balls that can be packed so that they all intersect
$B_0$, but are disjoint from each other. The value of $\kappa_1$ is
a lower bound on $\kappa_1^*$. In two dimensions $\kappa_1=5$ and
$\kappa_1^*=6$, though $\kappa_1=\kappa_1^*=12$ in three dimensions
and it seems likely the equality holds much more generally. In the
subsection, we assume $n$ is large enough so that \bea
\label{wide.enough.box} \sqrt{p}n^{1/p} > 2m \qmq{and} n^{1/p}>6,
\ena the first inequality being \eqref{n.large} specialized to the
unit radius case, and the second imposed so that
 Lemma \ref{lem:kappa} may be invoked over $C_n$.

For ${\cal A} \subset \{1,2,\ldots\}$ and ${\bf
d}=(d_\alpha)_{\alpha \in {\cal A}}$ a bounded collection of
positive integers, let $\kappa_{\bf d}$ be the maximum value of $k
\ge 0$ such that there exists $\Gamma \subset {\cal A}$ of size $k$
and a set of points ${\cal V}=(v_\alpha)_{\alpha \in \Gamma}$ in
$\mathbb{R}^p$ such that
\beas %\label{def:kappas}
B_\alpha \cap B_0 \not = \emptyset \qmq{and} \sum_{
    \gamma \in
\Gamma \setminus \{\alpha\} } {\bf 1}(B_\alpha \cap B_\gamma \not =
\emptyset) = d_\alpha-1 \qmq{for all $\alpha \in \Gamma$,} \enas
with $v_0$ the origin, and $B_\alpha$ the unit ball centered at
$v_\alpha$. That is, the constant $\kappa_{\bf d}$ is the maximum
size of a set $\Gamma$ of points in $\mathbb{R}^p$ such that the
number of neighbors of the points in the set $\Gamma$, among the
points indexed by $\Gamma \cup \{0\}$, drops from $d_\alpha$ to
$d_\alpha-1$ (increases from $d_\alpha-1$ to $d_\alpha$) upon the
removal (insertion) of the unit ball at the origin. If $d_\alpha=1$
for all $\alpha \in {\cal A}$, then $\kappa_{\bf d}=\kappa_1$.  Let
$d_{(\alpha)}$ be the values of $d_\alpha$ in a non-strict
decreasing order, that is, \beas d_{(1)} \ge d_{(2)} \ge \cdots
\enas

\begin{lemma} \label{lem:kappa}
In $\mathbb{R}^p$ for any dimension $p \ge 1$, \beas \kappa_{\bf d}
\le \sigma_{\bf d} \qmq{where} \sigma_{\bf d}= \sum_{\alpha \in
[\kappa_1]}  d_{(\alpha)},\enas and the bound is achieved when there
exists a pairwise disjoint collection of indices ${\cal Q}_\alpha
\subset {\cal A}$, $\alpha \in [\kappa_1]$, with $|{\cal
Q}_\alpha|=d_{(\alpha)}$ such that $d_\beta = d_{(\alpha)}$ for all
$\beta \in {\cal Q}_\alpha$.
\end{lemma}

\noindent {\em Proof:}  To show the upper bound, let ${\cal V}=
(v_\alpha)_{\alpha \in \Gamma}, \Gamma \subset {\cal A}$, be an
arbitrary collection of points in $\mathbb{R}^p$ such that, for all
$\alpha\in\Gamma$, the closed unit ball around $v_\alpha$ intersects
the closed unit ball around the origin and has exactly $d_\alpha-1$
neighbors  in ${\cal V}$. As the numbers $(d_\alpha)_{\alpha \in
{\cal A}}$ are bounded, the set ${\cal V}$ is finite. Let ${\cal R}$
be a subset of ${\cal V}$ of maximal size with the property that the
closed unit balls centered at the points of ${\cal R}$ are pairwise
disjoint. By the maximality of ${\cal R}$, each point in ${\cal V}
\setminus {\cal R}$ must be a neighbor of at least one point in
${\cal R}$, implying that ${\cal V}$ is the union of ${\cal R}$ and
those points of ${\cal V}$ that are a neighbor of some point of
${\cal R}$. Hence, \beas {\cal V} = \bigcup_{v \in {\cal R}} \left(
\{v\} \cup \bigcup_{w \in {\cal V}:\; w \sim v} \{w\} \right). \enas
As the point $v_\alpha \in {\cal V}$ has $d_\alpha-1$ neighbors in
${\cal V}$, and as $|{\cal R}|$ can be at most $\kappa_1$, we obtain
\beas |{\cal V}| \le \sum_{\alpha =1}^{\kappa_1} d_\alpha \le
\sum_{\alpha =1}^{\kappa_1} d_{(\alpha)} = \sigma_{\bf d}. \enas
Taking supremum over all such collections $(v_\alpha)_{\alpha \in
{\cal A}}$ we obtain the inequality $\kappa_{\bf d} \le \sigma_{\bf
d}$.

We now show that the bound is achieved when ${\bf d}$ satisfies the
given condition. By definition of $\kappa_1$, there exists a
collection of points $v_\alpha$, $\alpha \in [\kappa_1]$, in
$\mathbb{R}^p$ such that the closed unit balls $B_\alpha$, $\alpha
\in [\kappa_1]$ around each point intersect the closed unit ball
$B_0$ at the origin, but no other ball $B_\beta, \beta \in
[\kappa_1]\setminus \{\alpha\}$. Now consider the collection of
$\sigma_{\bf d}$ unit balls consisting of $d_{(\alpha)}$ copies of
the unit ball with center $u_\alpha$, for each $\alpha \in
[\kappa_1]$. Each of the $d_{(\alpha)}$ balls with center at
$u_\alpha$ has $d_{(\alpha)}$ neighbors when the closed unit ball at
the origin is included, but $d_{(\alpha)}-1$ neighbors when it is
not. Hence for such ${\bf d}$ we achieve $\kappa_{\bf d} \ge
\sum_{\alpha=1}^{\kappa_1} d_{(\alpha)}=\sigma_{\bf d}$. \bbox

\bigskip

We prove concentration for the neighborhood counts $Y_{ge}$
and $Y_{ne}$ by showing  that bounded size bias couplings for these variables exist with respective bounds
\begin{equation}\label{C.GG1}
        c_{ge} = |{\bf w}| |{\bf d}| \left( \sigma_{\bf d}+1 \right)
        \qmq{and} c_{ne}= |{\bf w}|\left(\sigma_{\bf d}+\sigma_{{\bf
                d}+1}+1\right),
        \end{equation}
where $|{\bf d}|$ and $|{\bf w}|$ are as in \eqref{def.absw.absd}
for given threshold and weight vectors, and where ${\bf d}+1$
denotes the vector $(d_{\alpha}+1)_{\alpha \in [m]}$.

Noting that the case $m=1$ is trivial, for $m=2$ we have $Y_{ge} \le 2|{\bf w}|$, hence $Y_{ge}^s$ is also
so upper bounded, and the inequality $Y_{ge}^s \le Y_{ge}+c_{ge}$
holds trivially, with similar remarks applying to $Y_{ne}$.  Hence
we may assume in the remainder of this section that $m \ge 3$. Under
\eqref{wide.enough.box}  it is not difficult to see that the
constant $\kappa_1$ computed over $C_n$ is the same as that over
$\mathbb{R}^p$, and Lemma \ref{lem:kappa} holds over $C_n$ as well.

For $\beta \not = \alpha$ and $x \in C_n$, letting \bea
\label{def:pbetau} p_\beta(x) = P(D(x,U_\beta) \le 2), \ena we have
that the conditional law ${\cal L}(M_\alpha|U_\alpha=x)$ is Poisson
Binomial ${\cal PB}({\bf p}_\alpha(u))$ where ${\bf
p}_\alpha(u)=(p_\beta(u))_{\beta \in [m] \setminus \{\alpha\}}$, and
$Y_{ge}$ and $Y_{ne}$ have expectations given respectively by \bea
\label{mu.GG1} \mu_{ge} = \sum_{\alpha \in [m]} w_\alpha \int_{C_n}
q_{ge}(d_\alpha,{\bf p}_\alpha(u))f_\alpha(u)du\,\, \mbox{and}\,\,
\mu_{ne} = \sum_{\alpha \in [m]} w_\alpha \int_{C_n}
q_{ne}(d_\alpha,{\bf p}_\alpha(u))f_\alpha(u)du. \ena

\begin{theorem} \label{thm:GGNN}
Concentration of measure inequalities
(\ref{a})-(\ref{bernsteincorollary}) hold for all $m \ge 1$ for $Y_{\rm ge}$ and $Y_{ne}$ as in \eqref{Y.occ}, computed on the neighbor count vector \eqref{Ma.GG1} for the germ-grain model, with corresponding $\mu$ and $c$ given respectively by \eqref{mu.GG1} and \eqref{C.GG1}.
\end{theorem}

\noindent {\em Proof:} First, we show that there exists a coupling of $Y_{ge}$ to $Y_{ge}^s$, having the
$Y_{ge}$-size biased distribution, that satisfies $Y_{ge}^s \le
Y_{ge} + c_{ge}$. To do so, for each $\alpha \in [m]$ we apply the
reasoning in the proof of Lemma~\ref{VolinRp}, replacing $m$ and
$U_0$ there by $m-1$ and $U_\alpha$. In particular, upon that
replacement \eqref{size.bias.mixture} becomes
\begin{multline} \label{size.bias.mixture.ngb}
{\cal L}({\cal U}|M_\alpha({\cal U}) \ge d_\alpha)= \int_{C_n} {\cal L}({\cal U}|M_{\alpha,u}({\cal U}) \ge d_\alpha) {\widetilde f}_\alpha (u)du, \\
\qmq{where} {\widetilde f}_\alpha(u)=\frac{P(M_{\alpha,u}({\cal U})
\ge d_\alpha)}{P(M_\alpha({\cal U}) \ge d_\alpha)} f_\alpha(u),
\end{multline}
where $M_{\alpha,u}({\cal U})$ is given by \eqref{Ma.GG1} with
$B_\alpha$ replaced by the unit ball centered at $u \in C_n$. Hence,
we sample ${\widetilde U}_\alpha$ on $C_n$ with density ${\widetilde
f}_\alpha(u)$, and for ${\widetilde U}_\alpha=u$, construct a
configuration ${\cal U}$ with law ${\cal L}({\cal
U}|P(M_{\alpha,u}({\cal U}) \ge d_\alpha)$ in order to achieve
${\cal L}({\cal U}|P(M_{\alpha}({\cal U}) \ge d_\alpha)$.

Continuing to follow the proof of  Lemma~\ref{VolinRp}, for each
$\alpha \in [m]$ and $a \in [m-1]_0$ we obtain coupled
configurations \beas {\cal U}_a^\alpha = \{U_{\beta,a}^\alpha, \beta
\in [m] \} \qmq{satisfying} {\cal L}({\cal U}_a^\alpha) = {\cal
L}({\cal U}|M_\alpha=a), \enas where  $U_{\alpha,a}^\alpha=u$ for
all $\alpha \in [m]$, and where ${\cal U}_a^\alpha$ and ${\cal
U}_{a+1}^\alpha$ differ in only one point indexed by, say,
$\beta_\alpha$, where $\beta_\alpha \not = \alpha$. We will say
$U_{{\beta}_a,a}^\alpha$ was removed from the configuration, into
which $U_{\beta_a,a+1}^\alpha$ is inserted.

Let $N_{\alpha,u}$ be constructed on this same space with
distribution ${\cal L}(M_{\alpha,u}({\cal U}))$. Applying
Corollary~\ref{cor:all.at.once} to $M_{u,\alpha}=M_{\alpha,u}({\cal
U})$, we obtain $A_{\alpha, u}$ satisfying \bea
\label{eq:bound.A.alpha} {\cal L}(N_{\alpha,u}+A_{\alpha,u}) = {\cal
L}(M_{\alpha,u}|M_{\alpha,u} \ge d_\alpha) \qmq{with} 0 \le
A_{\alpha,u} \le d_\alpha.\ena Recalling that $u$ was chosen with
density ${\widetilde f}_\alpha(u)$, and letting $N_\alpha$ and
$A_\alpha$ denote $N_{\alpha,u}$ and $A_{\alpha,u}$ respectively for
notational simplicity, by ${\cal L}(N_\alpha)={\cal
L}(M_{\alpha,u}({\cal U}))$ and \eqref{size.bias.mixture.ngb}, as in
Lemma~\ref{VolinRp} we have \bea \label{Ucount.right.conditional}
{\cal L}({\cal U}_{N_\alpha}^\alpha)={\cal L}({\cal U}) \qmq{and}
{\cal L}({\cal U}_{N_\alpha+A_\alpha}^\alpha)={\cal L}({\cal
U}|M_\alpha({\cal U}) \ge d_\alpha), \ena where ${\cal
U}_{N_\alpha}^\alpha$ and ${\cal
        U}_{N_\alpha+A_\alpha}^\alpha$ are configurations corresponding to
 $N_\alpha$ and $N_\alpha+A_\alpha$, respectively.

Let ${\bf N}_a^\alpha=(N_{\beta,a}^\alpha)_{\beta \in [m]}$ be the
occupancy counts corresponding to ${\cal U}_a^\alpha$ and write
$Y_{ge}({\cal V})$ for \eqref{Y.occ} evaluated on the configuration
${\cal V}$. By the second identity in
\eqref{Ucount.right.conditional}, the indicators $X_\beta^\alpha =
{\bf 1}(N_{\beta,N_\alpha+A_\alpha}^\alpha \ge d_\beta)$ satisfy
\eqref{Xbetaalphaiscond} with $X_\beta= {\bf
1}(N_{\beta,N_\alpha}^\alpha \ge d_\beta)$. Hence, by
Lemma~\ref{sblem}, with $I$ independent of all other variables with
distribution \eqref{def:Iprop}, $Y_{ge}^s=Y_{ge}({\cal
U}_{N_I+A_I}^I)$ has the $Y_{ge}=Y({\cal U}_{N_I}^I)$ size biased
distribution. The first identity in \eqref{Ucount.right.conditional}
shows that ${\cal L}({\cal U}_{N_I}^I)={\cal L}({\cal U})$.

As the removal of $U_{\beta_a,a}^\alpha$ from ${\cal U}_a^\alpha$
can only decrease the number of its neighbors in the configuration
${\cal U}_{a+1}^\alpha$, increases in occupancy counts for the
configuration ${\cal
    U}_{\beta_a,a+1}^\alpha$ over their values in ${\cal U}_a^\alpha$
can occur only for the point $U_{\beta_a,a+1}^\alpha$, taking the
place of $U_{\beta_a,a}^\alpha$, and its set of neighbors. However,
with $U_{\beta_a,a+1}^\alpha$ playing the role of the origin in
Lemma~\ref{lem:kappa}, its insertion into the point set
$\{U_{\beta,a}^\alpha: \beta \not = \beta_a\}$ can increase at most
$\sigma_{\bf d}$ of the counts $N_{\beta,a}^\alpha$, $\beta \not =
\beta_a$, of former value $d_\beta-1$ to counts
$N_{\beta,a+1}^\alpha$ of value $d_\beta$. Hence, also accounting
for the possible change of the count at the point indexed by
$\beta_a$, we have \beas |{\cal N}_{a+1}| \le  \sigma_{\bf d}+1,
\qmq{where} {\cal N}_{a+1} = \{\beta \in [m]: {\bf
1}(N_{\beta,a}^\alpha \ge d_\beta)=0,{\bf 1}(N_{\beta,a+1}^\alpha
\ge d_\beta)=1\}, \enas
 and therefore
\begin{multline*} \sum_{\beta \in [m]}w_\beta {\bf
1}(N_{\beta,a+1}^\alpha \ge d_\beta) = \sum_{\beta \not \in {\cal
N}_{a+1}}w_\beta {\bf
    1}(N_{\beta,a+1}^\alpha \ge d_\beta)+ \sum_{\beta \in {\cal N}_{a+1}}w_\beta {\bf
    1}(N_{\beta,a+1}^\alpha \ge d_\beta)\\
\le \sum_{\beta \not \in {\cal N}_{a+1}}w_\beta {\bf
1}(N_{\beta,a}^\alpha \ge d_\beta) + |{\bf w}||{\cal N}_{a+1}| \le
\sum_{\beta \in [m]}w_\beta {\bf
    1}(N_{\beta,a}^\alpha \ge d_\beta) + |{\bf w}|(\sigma_{\bf
    d} +1).
\end{multline*}
 Now,
by \eqref{eq:bound.A.alpha}, \beas Y_{ge}({\cal
U}_{N_\alpha+A_\alpha}^\alpha) \le Y({\cal U}_{N_\alpha}^\alpha) +
|{\bf w}| A_\alpha  ( \sigma_{\bf d}+1)  \le Y({\cal
U}_{N_\alpha}^\alpha) + |{\bf w}|  |{\bf d}| ( \sigma_{\bf d}
+1),\enas and mixing over $\alpha$ yields \beas
Y_{ge}^s=Y_{ge}({\cal U}_{N_I+A_I}^I) \le Y({\cal U}_{N_I}^I) +
|{\bf w}|  |{\bf d}| ( \sigma_{\bf d} +1) = Y_{ge} + |{\bf w}| |{\bf
d}| (\sigma_{\bf d}+1), \enas verifying the claim for $Y_{ge}$.

Next, we show that there exists a coupling of $Y_{ne}$ to $Y_{ne}^s$, having
the $Y_{ne}$-size biased distribution, satisfying $Y_{ne}^s \le
Y_{ne} + c_{ne}$, where $c_{ge}$ and $c_{ne}$ are given by
(\ref{C.GG1}). The construction for  $Y_{ne}$ will be similar, the only
difference being that the initial removal of $U_{\beta_a,a}^\alpha$
can cause $\sigma_{\bf d}$ counts to drop from $d_\alpha$ to
$d_\alpha-1$, while the insertion of $U_{\beta_a,a+1}^\alpha$ can
cause $\sigma_{{\bf d}+1}$ counts of value $d_\alpha$ to rise to
$d_\alpha+1$. Hence in this case we obtain, as claimed,
 \beas
 \sum_{\beta \in
[m]}w_\beta {\bf 1}(N_{\beta,a+1}^\alpha \ne d_\beta) \le
\sum_{\beta \in [m]}w_\beta {\bf 1}(N_{\beta,a}^\alpha \ne d_\beta)
+ |{\bf w}|(\sigma_{\bf d}+ \sigma_{{\bf d}+1}+1). \enas \bbox

\bigskip
When all points are uniformly distributed over $C_n$,
$f_\beta(u)=1/n$ and the probability~$p_\beta(u)$ in
\eqref{def:pbetau} is the constant $2^p\pi_p/n$ for all $u \in C_n$,
where $\pi_p$ is the volume of the unit ball in dimension $p$. Hence
with weights $w_\alpha=1$ for all $\alpha \in [m]$, we obtain \beas
\mu_{ge} = m P(M_\alpha \ge d) = m P(\mbox{Bin}(m-1,2^p\pi_p/n) \ge
d), \enas  with similar remarks applying to $\mu_{ne}$.

\subsection{Multinomial Occupancy} \label{mult:app} Among the many
applications of multinomial occupancy models, in which $n$ balls are
distributed independently to $m$ urns (see \cite{Kolchin78} for an
overview), are the well-known species trapping problem (see
\cite{Chao96}, \cite{Robbins68b}, or \cite{Starr79}) and the
closely-related problem of statistical linguistics (see
\cite{Efron76} and \cite{Thisted87}). The study of the number of
empty urns, or equivalently the $d=1$ case of $Y_{ge}$ in
\eqref{intro:YgeY=d}, was initiated in \cite{Re62}  and
\cite{Weiss58} where it was shown that the properly standardized
distribution of $Y_{ge}$ is asymptotically normal when balls land in
urns uniformly. Bounds in the $L^\infty$ metric between the standard
normal distribution and standardized finite sample distribution of
the $d=1$ case of $Y_{ge}$ was provided by \cite{Eng81} in the
uniform case, for $Y_{eq}$ by \cite{Pe09} in the uniform and some
non-uniform cases, and for all $d \ge 2$ for $Y_{eq}$ by
\cite{BaGo13} in the uniform case. Concentration of measure
inequalities for the number of empty urns were obtained in
\cite{dubhashi98}  by exploiting negative association, discussed in Section \ref{sec:NA}.

For $\alpha \in [m]$ let the component $M_\alpha$ of the vector
${\bf M}=(M_\alpha)_{\alpha \in [m]}$ count the number of balls in
urn $\alpha$ when $n$ balls are independently distributed into $m$
urns  and for $j \in [n]$ the location $L_j$ of ball $j$ is urn
$\alpha$ with probability $p_{\alpha,j}$. In particular,
\begin{equation}
\label{eq:counts.mult.occ} M_\alpha = \sum_{j=1}^n {\bf
1}(L_j=\alpha). 
\end{equation}

As in Section~\ref{sec:E-R}, we may assume that $p_{\alpha,j}<1$ for
all $(\alpha,j) \in [m] \times [n]$, that $\sum_{j=1}^n
p_{\alpha,j}>0$ for all urns $\alpha \in [m]$, and that each of the
summand indicators of $Y_{ge}$ and $Y_{ne}$ in \eqref{Y.occ} is
nontrivial. With ${\bf p}_\alpha=(p_{\alpha,j})_{j \in [n]}$ we have
$M_\alpha \sim {\cal PB}({\bf p}_\alpha)$ and, arguing as before,
the means $\mu_{ge}$ and $\mu_{ne}$ again have the form
\eqref{means.multinomial.occupancy}.

We may now summarize the main result of this subsection. 

\begin{theorem}\label{thm:multocc}
Concentration of measure inequalities
(\ref{a})-(\ref{bernsteincorollary}) hold for all $m \ge 1$ for multinomial occupancy counts,
\begin{enumerate}
  \item with
$Y_{ge}, \mu_{ge}$ and $c_{ge}$ given by \eqref{Y.occ},
\eqref{means.multinomial.occupancy}, and $|{\bf w}|$,
  \item with $Y_{ne}$, $\mu_{ne}$ and $c_{ne}$ given by \eqref{Y.occ},
\eqref{means.multinomial.occupancy} and $2|{\bf w}|$.
\end{enumerate}
\end{theorem}

Theorem~\ref{thm:multocc}  will follow immediately from Theorem \ref{thm:main} from the coupling construction provided by the next lemma.

\begin{lemma}  \label{lem:multinomial}
In the multinomial occupancy model there exists a coupling of $Y_{ge}$ to $Y_{ge}^s$, having the
$Y_{ge}$-size biased
 distribution, that satisfies $Y_{ge}^s \le Y_{ge} + |{\bf w}|$,
  and a coupling of $Y_{ne}$ to $Y_{ne}^s$, having the $Y_{ne}$-size biased distribution,
satisfying $Y_{ne}^s \le Y_{ne} + 2|{\bf w}|$.
\end{lemma}

\noindent {\em Proof:}  The reasoning of Lemma~\ref{lem:ergraph}
applies with only minimal changes. We take configurations in this
model to be \beas {\cal U}=\{L_j,  j \in [n]\}, \enas the collection
of locations of all $n$ balls. As for each $\alpha \in [m]$ the
corresponding count $M_\alpha$ has a LC marginal distribution with
support ${\cal S}_\alpha$ satisfying $\inf {\cal S}_\alpha=0$, in
order to invoke Theorem~\ref{+1-1:generally} it is only required to
show that, for all $\alpha \in [m]$, configurations $\{{\cal
U}_a^\alpha,a \in [m]\}$ exist with Properties $(0,\ge)$ and
$(1,\not =)$ satisfying \eqref{calU=dcalVa}.

As in Lemma~\ref{lem:ergraph}, through the use of
Lemma~\ref{add.subtract.one}, for each $\alpha \in [m]$ we obtain
configurations ${\cal U}_a^\alpha$ with corresponding counts ${\bf
N}_a^\alpha = (N_{\beta,a}^\alpha)_{\beta \in [m]}$ such that
$${\cal L}({\cal U}_a^\alpha)={\cal L}({\cal U}|M_\alpha=a) \qm{for all $a \in [n]_0$}
$$
and, for all $a \in [n-1]_0$, ${\cal U}_{a+1}^\alpha$ differs from
${\cal U}_a^\alpha$ by a single element, indexed by $j_a$, say. As
$N_{\alpha,a}^\alpha=a$ for all $a \in [n]_0$, we must have that
$L_{j_a,a}^\alpha \not =\alpha$ and $L_{j_a,a+1}^\alpha =\alpha$,
where $L_{j,a}^\alpha$ denotes the location of ball $j$ in
configuration ${\cal U}_a^\alpha$. In particular,
$N_{\beta,a+1}^\alpha \le N_{\beta,a}^\alpha$ for all $\beta \not =
\alpha$, and so our observation in \eqref{B.NA.case} guarantees that
the occupancy counts corresponding to $\{{\cal U}_a^\alpha,a \in
[m]\}$ have Property $(0, \geq)$.

For $Y_{ne}$, with $\{a,b\} \subset {\cal S}_\alpha$ satisfying
$|b-a|=1$, we note that the counts corresponding to the
configurations ${\cal U}_a^\alpha$ and ${\cal U}_b^\alpha$ differ
only at two indices, one of which is $\alpha$. Hence these counts
have Property $(1,\ne)$ by \eqref{bound.difference.by.B}. \bbox

\bigskip

In the asymptotic regime most studied, balls are uniformly
distributed, thresholds are constant and the weights are taken to be
identically 1. That is, $p_{\alpha, j} = 1/m, w_{\alpha} = 1$ and
$d_{\alpha} =d$ for each $\alpha \in [m]$ and $j \in [n]$.  For this
special case, the expectations in
\eqref{means.multinomial.occupancy} simplify to
\bea\label{mo:unifexpec} \mu_{ge} = m P(\mbox{Bin}(n,1/m) \geq d)
\quad \text{and} \quad \mu_{ne} = m(1 - P(\mbox{Bin}(n,1/m) = d)),
\ena and the concentration bounds obtained via size biasing can be used for $Y_{ge}$ with $c = 1$, and
for $Y_{ne}$ with $c=2$.

The expectations in (\ref{mo:unifexpec}) are of a form similar to
those of (\ref{erstandardmeans}) for the standard
Erd\H{o}s-R\'{e}nyi model. Thus, by arguing as in Section
\ref{sec:E-R} we can study in the same manner the behavior of the
bounds obtained here.

\begin{remark}
(Multivariate hypergeometric sampling ) 
The techniques used here to construct a bounded coupling can be modified to obtain concentration inequalities for another relevant model where we have negative association. Let $n$ be the sum of the given positive
integers $(n_\alpha)_{\alpha \in [m]}$, and consider an urn
containing $n$ colored balls, $n_\alpha$ of which are of color
$\alpha$. For $s\in[n]_0$ let $M_\alpha$ be the number of balls of
color $\alpha$ obtained upon sampling $s$ distinct balls uniformly
from the urn without replacement, and set ${\bf
M}=(M_\alpha)_{\alpha \in [m]}$. Let  $Y_{ge}$ and $Y_{ne}$ be as in
\eqref{Y.occ}.

Then, for all $\alpha \in [m]$ the distribution of
 $M_\alpha$ is hypergeometric, and the
expected values of $Y_{ge}$ and $Y_{ne}$ are, respectively,
\begin{equation}\label{mu.hyper}
  \mu_{ge}  = \sum_{\alpha \in [m]} w_{\alpha} \sum_{j\ge d_{\alpha}}q(j;n_\alpha,s,n) \qmq{and}
    \mu_{ne}= \sum_{\alpha \in
[m]} w_{\alpha} \sum_{j \neq d_{\alpha}} q(j;n_\alpha,s,n)
\end{equation} where
\begin{equation}\label{dhyper} q(j;k,\ell,i)= \binom{k}{j}
\binom{i-k}{\ell-j}\left/\binom{i}{\ell}\right.
\end{equation}
for values of $j,k,\ell,i$ such that the quotient \eqref{dhyper} is defined, setting $q(j;k,\ell,i)=0$ otherwise. Moreover, 
 it is shown in
\cite[Theorem~A]{Ehm91} that a hypergeometric random variable can be
written as a sum of independent but non-identically-distributed
Bernoulli random variables, from which we may conclude  that the hypergeometric distribution is LC.  Our techniques above then can be modified to show that
there exists a coupling of $Y_{ge}$ to $Y_{ge}^s$, having the
$Y_{ge}$-size biased distribution, that satisfies $Y_{ge}^s \le
Y_{ge} + |{\bf w}|$, and a coupling of $Y_{ne}$ to $Y_{ne}^s$,
having the $Y_{ne}$-size biased distribution, satisfying $Y_{ne}^s
\le Y_{ne} +2 |{\bf w}|$. These couplings provide concentration of measure inequalities where the means  are  given by \eqref{mu.hyper}. We omit the details.
\end{remark}

\section{Comparisons} \label{TalAzu}
In this section, we compare our results to concentration bounds obtained by other
means. Our comparisons will be with the following three well known
techniques: (i) McDiarmid's Inequality, (ii) Use of negative
association and (iii) Self Bounding and Certifiable functions. Of
these three, the last technique is the most comparable. For simplicity and concreteness, in
most of our comparisons below we
will consider 
\begin{equation}
\label{unitwconstd} Y_{ge} = \sum_{\alpha\in [m]}
\mathbf{1}(M_{\alpha} \geq d)
\end{equation}
with unit weighting and constant threshold count.

\subsection{McDiarmid's Inequality}\label{sec:McDiarmid}
One of the most useful concentration results is the McDiarmid, or
bounded difference, inequality which is a consequence of the
Azuma-Hoeffding bound; see \cite{hoe}, \cite{azu} and \cite{mcd}.
The inequality applies to quantities $Y$ that can be expressed as a
function $f(X_1,\ldots,X_n)$ of independent random variables
$X_1,\ldots,X_n$ when, for all $i \in [n]$, there exists a constant
$c_i$ such that
 \bea \label{ci:AH}
\sup_{x_i,x_i'}|f(x_1,\ldots,x_i,\ldots,x_n) -
f(x_1,\ldots,x_i',\ldots,x_n)| \le c_i. \ena

Under these conditions, the inequality provides the right tail bound
\bea \label{AH:ineq} P(Y - E Y \ge t) \le \exp\left(-\frac{2t^2}{
\sum_{i=1}^n c_i^2}\right), \ena and a corresponding left tail
bound.

Although the bounded difference inequality is powerful and easy to
apply, the quantity $\sum_{i=1}^n c_i^2$ on which it depends,
obtained by taking supremums in \eqref{ci:AH} to estimate the
worst-case behavior of $f$, may not accurately reflect the
concentration properties of $f$.

To take the simplest example, let $Y$ have the Binomial distribution
$\mbox{Bin}(n,p)$. As $Y$ can be written as the sum of independent
Bernoullis, inequality \eqref{ci:AH} is satisfied with $c_i=1$ and
\eqref{AH:ineq} yields \bea \label{AH.n} P(Y - np \ge t) \le
\exp\left(-\frac{2t^2}{n}\right). \ena However, for the Binomial it
is known (see \cite{mcdiarmid2006concentration}, for instance) that
the true decay rate is $\exp(-t^2/2np)$. In particular, use of
\eqref{AH.n} may not be adequate in situations where $p$ is small.

Applying Lemma~\ref{sblem} to $Y$, represented as an independent sum
of indicators, we find that $Y^s$ can be formed by replacing any of
the summand indicators by 1, yielding $Y^s \le Y+1$. Hence, the bound
\eqref{bernsteincorollary} yields
$$P\left(Y-np \ge t\right)\le \exp\left(-\frac{t^2}{2(np+t/3)}\right)
\qm{for all $t>0$,}$$ which, specializing to the case $p \in
(0,1/4)$, improves on the Azuma-Hoeffding bound
 (\ref{AH.n}) in the range $0<t<3n(1/4-p)$, with upper range increasing to $(0,\infty)$ as $n \rightarrow \infty$.

We now turn  to the standard Erd\H{o}s-R\'{e}nyi random graph
$\mathcal{G}_m$ on $m$ vertices with fixed edge probabilities $p$, as
considered in Section~\ref{sec:E-R}, and let $Y_{ge}$ be given by
\eqref{unitwconstd} where $M_{\alpha}$ is the degree of vertex
$\alpha$. Clearly $Y_{ge}$ can be written as a function~$f$ of
$n={\binom{m}{2}}$ independent indicators $X_1,\ldots,X_n$, where
$X_i$ denotes the presence of a given edge with respect to some
fixed labeling. As a change in any  $X_i$ affects the degree of
exactly two vertices, $f$ satisfies the bounded differences
condition (\ref{ci:AH}) with $c_i=2$ for each $i =1,\ldots, n$.
Hence \eqref{AH:ineq} and the complementary left tail inequality
yield \bea\label{AHERnonis} \max\left\{P(Y_{ge}-\mu_{ge} \leq -t),
P(Y_{ge}-\mu_{ge} \geq t) \right\} \leq
\exp\left(-\frac{t^2}{m(m-1)}\right), \ena where $\mu_{ge}$  is
given by \eqref{erstandardmeans}. Comparing the left tail bounds of
\eqref{AHERnonis} with \eqref{sbERnonis}, we see the size bias bound
is preferred when \beas \mu_{ge} \le m(m-1)/(2d+2). \enas As
$\mu_{ge} \le m$ due to $Y$ being the sum of
$m$ indicators, this inequality is always satisfied for $m \ge
2d+3$, and  we see that the order of the exponent is improved from
$O(-t^2/m^2)$ to $O(-t^2/m)$. Similar improvements will also hold
for the Erd\H{o}s-R\'{e}nyi type graph models with inhomogeneous
edge probabilities, which, depending on their values, can become
even more significant. 

\subsection{Negative Association}\label{sec:NA} Negative association has been used
successfully to obtain concentration of measure inequalities for
occupancy models. We recall from \cite{Joag-Dev83}  (see
also \cite{dubhashi98} and \cite{Shaked}) that a family of random variables
$X_1,X_2,\ldots,X_m$ is said to be negatively associated if for any
disjoint subsets $A_1,A_2 \subset [m]$,
$$ E(f(X_i; i \in A_1) g(X_j ; j \in A_2)) \leq  E(f(X_i; i \in A_1)) E(g(X_j; j \in A_2))$$
whenever $f$ and $g$ are coordinate-wise nondecreasing functions for
which these expectations exist.

For the multinomial occupancy model of Section~\ref{mult:app},
the results of \cite{Joag-Dev83} show that the indicator
summands ${\bf 1}(M_\alpha \ge d_\alpha)$ of $Y_{ge}$ are negatively
associated. Referring to Proposition 7 of
\cite{dubhashi98}, when $X_1,X_2,\ldots,X_m$ are negatively
associated indicators, the random variable $Y = \sum_{i=1}^m X_i$
satisfies the right tail bound of \eqref{arba:bound}, and hence
bounds \eqref{a} and \eqref{bernsteincorollary} with $c=1$.  Thus,
for this case,  the estimates obtained for the right tail via
negative association are at least as good as the ones that are
obtained by using size biasing, with the same holding for the left
tail. Indeed, Chernoff's bound for sums of independent random
variables (\cite{Bo13}, page 24) remains true for negatively
associated sums, indicating that one may have strict improvements
over the size bias method.

However, none of the other statistics discussed here can be handled using negative association. For instance, one
cannot use negative association for our applications to random
graphs and germ-grain models in Sections~\ref{sec:E-R} and
\ref{sec:GG}. In particular, for the standard Erd\H{o}s-R\'{e}nyi random graph, a simple application of Harris'
inequality shows that the summand variables of $Y_{ge}$ are
positively associated.

Moreover, even for the multinomial occupancy model where negative
association can be used for $Y_{ge}$, the indicator summands in the
multinomial occupancy count
\beas %\label{neqcase}
Y_{ne} = \sum_{\alpha \in [m]} \mathbf{1}(M_{\alpha} \neq
d_{\alpha}) \enas are not negatively associated when the thresholds
$d_{\alpha}$ are not all $0$. Hence, the method of \cite{dubhashi98}
no longer applies, while the methods in this paper are still valid. For instance, when
$d_{\alpha} =1$ for each $\alpha \in [m]$ and balls are distributed
uniformly, Part 2 of Theorem \ref{thm:multocc} yields that \eqref{arba:bound}
holds with $c=2$ and
$$\mu_{ne}=  m(1-P(M_1 =1))=m\left(1-\frac{n}{m}\left(1-\frac{1}{m}\right)^{n-1}\right).
 $$

\subsection{Self-bounding and Certifiable Functions}\label{sec:self.bd} We have seen above
that bounds produced by size biasing may improve on the bound
\eqref{AH:ineq} obtained using the bounded difference inequality as
it replaces the sum $\sum_{i=1}^n c_i^2$ by some function of the
mean of $Y$. Bounds produced by the method of self bounding
functions \cite{mcdiarmid2006concentration}, of which certifiable
functions are a special case, also have this advantage. We focus on
the latter, as it is more straightforward to address the
applications studied here in the framework of certifiable functions.

We begin by recalling the relevant definitions and results on
certifiable functions from \cite{mcdiarmid2006concentration}; see
also \cite{Bo13}. Let $c>0$, $a \geq 0$, and $b$ be given, and let a
nonnegative measurable function $f$ on the product space $\Omega =
\Pi_{i =1}^n\Omega_i$ satisfy the following two conditions.
\begin{enumerate}

\item[(i)] For each $x \in \Omega$, changing
any coordinate $x_j$ changes the value of $f(x)$ by at most $c$.

\item[(ii)] If $f(x) =s$ then there is a set of coordinates $C \subset [n]$ of size at most $as +b$
that \textit{certifies} $f(x) \geq s$. That is, if the coordinates
$i \in C$ of $y \in \Omega$ agree with those of $x$, then $f(y) \ge
s$.

\end{enumerate}

Let $X_1,\ldots,X_n$ be independent random variables with $X_i$
taking values in $\Omega_i$, $Y=f(X_1,...,X_n)$ where $f$ satisfies
(i) and (ii) above, and $\mu=EY$. Then for all $t \geq 0,$
\begin{multline}\label{certboth}
P(Y - \mu \leq -t) \leq \exp\left(-\frac{t^2}{2c^2 (a
\mu + b + t / 3c)}\right) \quad \mbox{and}\\
    P(Y - \mu \geq t)\leq \exp\left(- \frac{t^2}{2c^2(a \mu + b+
    at)}\right).
\end{multline}

Before moving to a discussion of specific examples, we note that the
asymptotic Poisson order $O(\exp(-t \log t))$ as $t \rightarrow
\infty$ of the bound \eqref{arba:bound} with $c=1$ and $\mu=1$, is
superior to the order $O(\exp(-t))$ of the bound \eqref{certboth}
with $c=1$ and $a=1/2$, say, with similar types of improvement in
order holding for other choices of constants. The order of the
bounds achieved by certifiable functions, and self bounding
functions more generally, seems to be intrinsic. Regarding using the
entropy method to prove concentration inequalities for self bounding
functions, via log Sobolev inequalities in particular, the authors
of \cite{Bo13} note after the proof of Theorem 6.21 that,  `At least
for $a>1$, there is no hope to derive Poissonian bounds\ldots  for
the upper tail.'

\begin{comment}
To focus on a specific example, consider the multinomial occupancy
model of Section~\ref{mult:app}, and let $Y_{ge}$ be given by
\eqref{unitwconstd} where $M_{\alpha}$ is the number of balls in urn
$\alpha$. The variable $Y_{ge}$ can clearly be written as a function
$f$ of the  independent locations $L_j$ of ball  $j \in [n]$. It is
not difficult to verify that $f$ is certifiable with $c=1, a=d$ and
$b=0$. Thus, from \eqref{certboth},
 $Y_{ge}$ satisfies
\begin{equation}\label{certexampleboth}
P(Y_{ge}-\mu_{ge} \leq -t) \leq \exp
\left(-\frac{t^2}{2(d\mu_{ge}+t/3)}\right)\quad \mbox{and}\\
P(Y_{ge}-\mu_{ge} \geq t) \leq \exp
\left(-\frac{t^2}{2d(\mu_{ge}+t)}\right).
\end{equation}
Applying size biasing, Part~3 of Theorem~\ref{thm:monster2} shows
that the lower and upper tail bounds \eqref{a} and
\eqref{bernsteincorollary} hold for $Y_{ge}$ with $c=1$, and a
simple computation shows that both these inequalities strictly
outperform their counterparts in (\ref{certexampleboth}).

Similar remarks apply to the statistic $Y_{ne} =
 \sum_{\alpha \in [m]} \mathbf{1}(M_{\alpha} \neq d)$, which is a certifiable
function with $c=2, b=n$ and $a=0$. The left tail size bias bound
$\exp(-t^2/(4\mu_{ne}))$ obtained via (\ref{a}), will always be
preferable to the left tail bound obtained from \eqref{certboth}, as
$4\mu_{ne} \le 4n \le 8(n+t/6)$.
\end{comment}

To focus on a specific example, consider the random graph model of
Section~\ref{sec:E-R}, and let $Y_{ge}$ be given by
\eqref{unitwconstd} where $M_{\alpha}$ is the degree of vertex
$\alpha$ and $d \ge 2$. One can now easily show that the statistic
$Y_{ge}$ is certifiable with
 $c=2, a =d$ and $b =0$, though both the lower and upper tail
 bounds \eqref{a} and \eqref{bernsteincorollary}
are superior to those obtained via \eqref{certboth}.

Finally, we note that the Poisson tail concentration of measure
inequalities of Theorem~\ref{thm:main} will always provide further
improvements over the bounds (\ref{a}) and
(\ref{bernsteincorollary}) applied in the previous paragraphs.
However, the form of these latter bounds, being simpler than that of
\eqref{arba:bound},  allow for an  easier comparison  with
(\ref{certboth}), and although  they are not the strongest bounds of
those produced by the size bias method, they still suffice to
demonstrate the improvements claimed.

\appendix
\section{Monotonicity Assumption } \label{proof:main}
We prove that the monotonicity assumption that $Y^s \ge Y$ assumed
in \cite{GhGo11a} and \cite{GhGo11b} for the left tail bound can be
removed, and that only $Y^s \le Y+c$ is required for (\ref{a}).
First, we may assume that $Y$ is not almost surely constant as
inequality (\ref{a}) is trivially satisfied in that case. Since $Y
\ge 0$ a.s., for all $\theta <0$ the moment generating function
$m(\theta)=E(e^{\theta Y})$ of $Y$ exists in an open interval
containing $\theta$ and is differentiable at $\theta$.
Differentiating under the expectation by dominated convergence and
then applying the characterization of the size bias distribution
(\ref{def:sb}), followed by an application of the inequality $1+x
\leq e^x $, we obtain
\begin{multline} \label{monremove1}
m'(\theta) = E(Ye^{\theta Y}) = \mu  E(e^{\theta Y^s}) = \mu
E(e^{\theta Y} e^{\theta (Y^s-Y)} ) \geq  \\
\mu E(e^{\theta Y} (1+\theta (Y^s-Y))) \ge \mu E(e^{\theta Y}
(1+\theta c )) = \mu  (1+\theta c) m(\theta),
\end{multline}
where we have used $Y^s-Y\leq c$ and $\theta <0$.  Rearranging terms
in (\ref{monremove1}) yields \bea\label{monremove2} 0 \leq
m'(\theta) - \mu (1+\theta c) m(\theta), \ena and multiplying each
side of (\ref{monremove2}) by $e^{-\mu (\theta + c \theta^2/2)}$ we
see that
\begin{eqnarray}\label{monremove3}
0 \leq
%m'(\theta)e^{-\mu (\theta +c \theta^2/2)}-\mu m(\theta) (1+\theta c) e^{-\mu (\theta +c \theta^2/2)} =
(m(\theta) e^{-\mu (\theta +c \theta^2/2)})' \qmq{for all
$\theta<0$.}
\end{eqnarray}
Integrating both sides of (\ref{monremove3}), and using $m(0)=1$,
yields
$$0 \leq \int_{\theta}^0 (m(x) e^{-\mu
(x +c x^2/2)})'dx = 1 - m(\theta) e^{-\mu (\theta + c \theta^2
/2)}$$ and hence \bea\label{monremove4}m(\theta) \leq e^{\mu (\theta
+c \theta^2/2)}.\ena

Next letting $M(\theta)= E(e^{\theta(Y-\mu)})=e^{- \mu \theta } m
(\theta)$ and applying (\ref{monremove4}), we obtain the bound
$$M(\theta) \leq e^{-\mu \theta} e^{\mu(\theta +c \theta^2/2)}=
e^{\mu c \theta^2 /2}.$$ Hence for fixed $t >0$ and all $\theta <0$,
$$P(Y-\mu  \leq -t) = P(e^{\theta (Y-\mu)}\geq e^{-\theta t}) \leq e^{\theta t} M(\theta) \leq e^{\theta t +  \mu c \theta^2 /2}$$ by Markov's
inequality. Substituting $\theta = -t/c \mu$ yields  inequality
\eqref{a}.\bbox

\end{document}